# Fréchet and Mordukhovich Derivative (Coderivative) and Covering Constant for Single-Valued Mapping in Euclidean Space with Application (II)

Jinlu Li

Department of Mathematics
Shawnee State University
Portsmouth, OH 45662 USA
Email: jli@shawnee.edu

**Abstract** We continue the study in [15] for calculating the Fréchet derivatives and Mordukhovich derivatives (coderivatives) and covering constants for single-valued mappings in Euclidean spaces ([15] is part I). In this paper, we particularly consider a norm preserving mapping $f: \mathbb{R}^2 \to \mathbb{R}^2$ that is defined by (1.1) in Section 1. We will find the precise solutions of Fréchet derivative and Mordukhovich derivative at every point in $\mathbb{R}^2$. By using these solutions, we will find the covering constant for this mapping $f$ is exact 1 at every point in $\mathbb{R}^2$ except the origin. Then we extend this mapping to $\mathbb{R}^4$. Finally, by using the covering constant for $f$ and by applying the Arutyunov Mordukhovich and Zhukovskiy Parameterized Coincidence Point Theorem, we will solve some parameterized equations.

**Keywords**  Fréchet derivative; Mordukhovich derivative (or coderivatives); covering constant; AMZ Theorem

## 1. Introduction

The Mordukhovich derivative (or the Mordukhovich coderivative) for set-valued mappings lays the foundation and it plays fundamental and crucial role in the theory of generalized differentational analysis in Banach spaces (See [15−20]). The theory of generalized differentiation in set-valued analysis has been rapidly developed during the last 30 years and has been applied to many fields (see [1−8, 11−14]). In contrast to the differentiation theory in calculus, in some real applications involved with some considered mappings, it is required to know the exact solution of its Mordukhovich derivative. Hence, our concern is that for an exactly given mapping between Banach spaces, how do we find its precise solutions of the Mordukhovich derivative, which may be considered as some Mordukhovich derivative formulas? For this purpose, in [9−10] the present author studied the standard metric projection operator in Hilber spaces and uniformly convex and uniformly smooth Banach spaces and find the solutions of its Mordukhovich derivative to some special closed and convex subsets (including closed balls).

In this paper, we consider a norm preserving mapping $f$ in $\mathbb{R}^2$, which is defined by (1.1) below. We will find its precise solutions of the Mordukhovich derivative, by which, we will find the exact covering constant for this mapping $f$.

Let $(\mathbb{R}^2, \|\cdot\|)$ be the standard 2-d Euclidean space with the ordinal Hilbert $L_2$-norm and row vectors. Let $\theta$ denote the origin of $\mathbb{R}^2$. In this paper, we consider a single-valued mapping $f: \mathbb{R}^2 \to \mathbb{R}^2$, which is studied in Examples 2 in [4] and 4.2 in [2]. $f$ is defined by

$$f((x_1, x_2)) = \left( \frac{x_1^2 - x_2^2}{\sqrt{x_1^2 + x_2^2}}, \frac{2x_1 x_2}{\sqrt{x_1^2 + x_2^2}} \right), \text{ for } (x_1, x_2) \in \mathbb{R}^2 \setminus \{\theta\} \text{ with } f(\theta) = \theta. \tag{1.1}$$

This mapping has many interesting properties that will be investigated in Section 3. One of the properties $f$ is norm preserving. That is that, for any $x \in \mathbb{R}^2$, one has $\|f(x)\| = \|x\|$.

For $\bar{x}, \bar{y} \in \mathbb{R}^2$ with $\bar{y} = f(\bar{x})$, let $\text{cov} f(\bar{x}, \bar{y})$ denote the exact covering bound of $f$ around $(\bar{x}, \bar{y})$. Let $\alpha(f, \bar{x}, \bar{y})$ denote the supremum of locally $\alpha-$ covering of $f$ in a neighborhood of $(\bar{x}, \bar{y})$. Let $\hat{\alpha}(f, \bar{x}, \bar{y})$ denote the covering constant for $f$ at point $(\bar{x}, \bar{y})$. All these concepts and some of their properties will be reviewed in Section 2. The concepts of the exact covering bound and covering constant for mappings have been widely applied to solve some inverse problems (see [4]), to solve some coincidence point problems (see [1]) and some other problems (see [11−14]). Hence, to find (to calculate) the exact covering bound and covering constant for the considered mappings in some application problems becomes an extremally important object in variational analysis, optimization theory and related fields.

However, notice that (see Section 2) the concepts of the exact covering bound $\text{cov} f(\bar{x}, \bar{y})$ and the supremum $\alpha(f, \bar{x}, \bar{y})$ locally $\alpha-$ covering of $f$ around $(\bar{x}, \bar{y})$ are geometrically defined by the graphs of the considered mapping $f$. So, they may be geometrically calculated based on the graph of $f$. In [4], the authors gave an elegant (geometrical) proof for the following result regarding to the supremum of locally $\alpha-$ covering of $f$:

$$\alpha(f, \bar{x}, \bar{y}) = 1, \text{ for } \bar{x}, \bar{y} \in \mathbb{R}^2 \text{ with } \bar{y} = f(\bar{x}). \tag{1.2}$$

The covering constant for $f$ at point $(\bar{x}, \bar{y})$, $\hat{\alpha}(f, \bar{x}, \bar{y})$, is defined by Mordukhovich derivatives of $f$ at point around $(\bar{x}, \bar{y})$. It is well-known that it is very difficult to precisely find the Mordukhovich derivatives of a single-valued mapping at a point, which deduce the difficulty for finding the covering constant $\hat{\alpha}(f, \bar{x}, \bar{y})$ for $f$ at point $(\bar{x}, \bar{y})$.

It is gratifying that the close connection (equality property) between $\text{cov} f(\bar{x}, \bar{y})$, $\alpha(f, \bar{x}, \bar{y})$ and $\hat{\alpha}(f, \bar{x}, \bar{y})$ in Theorem 4.1, Corollaries 4.2, 4.3 and Proposition 3.6 in [22], from which we have that

$$\hat{\alpha}(f, \bar{x}, \bar{y}) = \text{cov } f(\bar{x}, \bar{y}) = \alpha(f, \bar{x}, \bar{y}) = 1, \text{ for } \bar{x}, \bar{y} \in \mathbb{R}^2 \text{ with } \bar{y} = f(\bar{x}). \tag{1.3}$$

In this paper, we use the results in [15] to precisely find the Fréchet derivatives and Mordukhovich derivatives of $f$ at every point in $\mathbb{R}^2$ and to prove that the Mordukhovich derivatives of $f$ are norm preserving. Then, by these properties, we prove that the covering constant $\hat{\alpha}(f, \bar{x}, \bar{y}) = 1$, for every $\bar{x}, \bar{y} \in \mathbb{R}^2$ with $\bar{y} = f(\bar{x})$.

The concept of covering constant plays an important and crucial role in Arutyunov Mordukhovich and Zhukovskiy Parameterized Coincidence Point Theorem (It is simply named as AMZ Theorem), which is considered as one of the most important applications of covering constants. In Section 5, we will use the result that $\hat{\alpha}(f, \bar{x}, \bar{y}) = 1$ to solve some parameterized equations.

## 2. Preliminaries

Let $(X, \|\cdot\|_X)$ and $(Y, \|\cdot\|_Y)$ be real Banach spaces, which have topological dual spaces $(X^*, \|\cdot\|_{X^*})$ and $(Y^*, \|\cdot\|_{Y^*})$, respectively. Let $\langle \cdot, \cdot \rangle_X$ denote the real canonical pairing between $X^*$ and $X$ and $\langle \cdot, \cdot \rangle_Y$ the real canonical pairing between $Y^*$ and $Y$. Let $\theta_X$ and $\theta_Y$ denote the origins in $X$ and $Y$, respectively. For any $x \in X$ and $r > 0$, let $B_X(x, r)$ and $S_X(x, r)$ denote the closed ball and sphere in $X$ centered at point $x$ with radius $r$, respectively.

Let $g: X \to Y$ be a single-valued mapping and let $\bar{x} \in X$. If there is a linear and continuous mapping $\nabla g(\bar{x}): X \to Y$ such that

$$\lim_{x \to \bar{x}} \frac{g(x) - g(\bar{x}) - \nabla g(\bar{x})(x - \bar{x})}{\|x - \bar{x}\|_X} = \theta_Y,$$

then $g$ is said to be Fréchet differentiable at $\bar{x}$ and $\nabla g(\bar{x})$ is called the Fréchet derivative of $g$ at $\bar{x}$.

Since this paper only deals with the single-valued mapping $f$ defined in (1.1), so we only review the concepts of Mordukhovich derivative (or coderivative) for single-valued mappings in Banach spaces (see [15–20] for more details). Let $\Delta$ be a nonempty subset of $X$ and let $g: \Delta \to Y$ be a single-valued mapping. For $x \in \Delta$ and $y = g(x)$, a set-valued mapping $\widehat{D}^* g(x, y): Y^* \to X^*$ is defined, for any $y^* \in Y^*$, by (see Definitions 1.13 and 1.32 in Chapter 1 in [15])

$$\widehat{D}^* g(x, y)(y^*) = \left\{ z^* \in X^*: \limsup_{\substack{(u, g(u)) \to (x, g(x)) \\ u \in \Delta}} \frac{\langle z^*, u - x \rangle_X - \langle y^*, g(u) - g(x) \rangle_Y}{\|u - x\|_X + \|g(u) - g(x)\|_Y} \leq 0 \right\}.$$

If

$$\widehat{D}^* g(x, y)(y^*) \neq \emptyset, \text{ for every } y^* \in Y^*, \tag{2.1}$$

then $g$ is said to be Mordukhovich differentiable (co differentiable) at the point $x$ and $\widehat{D}^* g(x, y)$ is called the Mordukhovich derivative (which is also called Mordukhovich coderivative, or coderivative) of $g$ at point $x$. For this single-valued mapping $g: \Delta \to Y$, we write

$$\widehat{D}^* g(x, y)(y^*) \equiv \widehat{D}^* g(x)(y^*).$$

Furthermore, if $g: \Delta \to Y$ is a single-valued continuous mapping. Then, by the above definition, the Mordukhovich derivative of $g$ at point $x$ is calculated by

$$\widehat{D}^* g(x)(y^*) = \left\{ z^* \in X^*: \limsup_{\substack{u \to x \\ u \in \Delta}} \frac{\langle z^*, u - x \rangle_X - \langle y^*, g(u) - g(x) \rangle_Y}{\|u - x\|_X + \|g(u) - g(x)\|_Y} \leq 0 \right\}, \text{ for any } y^* \in Y^*. \tag{2.2}$$

The following theorem shows the connection between Fréchet derivatives and Mordukhovich derivatives for sing-valued mappings. The results of the following theorem provide a powerful tool to calculate the Mordukhovich derivatives by the Fréchet derivatives of single-valued mappings.

**Theorem 1.38 in [15]**. *Let $X$ be a Banach space with dual space $X^*$ and let $g: X \to Y$ be a single-valued mapping. Suppose that $g$ is Fréchet differentiable at $x \in X$ with $y = g(x)$. Then, the Mordukhovich derivative of $g$ at $x$ satisfies the following equation*

$$\widehat{D}^* g(x)(y^*) = \{(\nabla g(x))^*(y^*)\}, \text{ for all } y^* \in Y^*.$$

One of the important applications of Mordukhovich derivatives of set-valued mappings is to define the covering constants for set-valued mappings. The covering constant for $\Phi: X \rightrightarrows Y$ at point $(\bar{x}, \bar{y}) \in \text{gph } \Phi$ is defined by (see (2.4) in [1])

$$\hat{a}(\Phi, \bar{x}, \bar{y}) := \sup_{\eta > 0} \inf\{\|z^*\|_{X^*}: z^* \in \widehat{D}^* \Phi(x, y)(w^*), x \in \mathbb{B}_X(\bar{x}, \eta), y \in \Phi(x) \cap \mathbb{B}_Y(\bar{y}, \eta), \|w^*\|_{Y^*} = 1\}. \tag{2.3}$$

Here, $\|\cdot\|_{X^*}$ and $\|\cdot\|_{Y^*}$ denote the norms in $X^*$ and $Y^*$, respectively. $\mathbb{B}_X(\bar{x}, \eta)$ is the closed ball in $X$ centered at $\bar{x}$ with radius $\eta$, and $\mathbb{B}_Y(\bar{y}, \eta)$ is the closed ball in $Y$ centered at $\bar{y}$ with radius $\eta$.

In particular, let $g: X \to Y$ be a single-valued mapping. For any $\bar{x}, \bar{y} \in X$ with $\bar{y} = g(\bar{x})$, (2.4) becomes to

$$\hat{\alpha}(g, \bar{x}, \bar{y}) = \sup_{\eta > 0} \inf\{\|z^*\|_{X^*} : z^* \in \widehat{D}^* g(x)(w^*), x \in \mathbb{B}_X(\bar{x}, \eta), g(x) \in \mathbb{B}_Y(\bar{y}, \eta), \|w^*\|_{Y^*} = 1\}. \quad (2.4)$$

In this paper, we consider Hilbert spaces as special cases of Banach spaces. Let $(H, \|\cdot\|)$ be a real Hilbert space with inner product $\langle \cdot, \cdot \rangle$ and origin $\theta$. Let $g: H \to H$ be a single-valued mapping and let $\bar{x} \in X$. Then, the Fréchet derivative of $g$ at $\bar{x}$ is defined by $\nabla g(\bar{x}): H \to H$ such that

$$\lim_{x \to \bar{x}} \frac{g(x) - g(\bar{x}) - \nabla g(\bar{x})(x - \bar{x})}{\|x - \bar{x}\|} = \theta.$$

Let $\Delta$ be a nonempty subset of $H$ and let $g: \Delta \to H$ be a single-valued continuous mapping. For $x \in \Delta$ and $y = g(x)$, by (2.2), the Mordukhovich derivative of $g$ at point $x$ becomes

$$\widehat{D}^* g(x)(y) = \left\{ z \in H : \limsup_{\substack{u \to x \\ u \in \Delta}} \frac{\langle z, u - x \rangle - \langle y, g(u) - g(x) \rangle}{\|u - x\| + \|g(u) - g(x)\|} \leq 0 \right\}, \text{ for any } y \in H. \quad (2.5)$$

In this case, for any $\bar{x}, \bar{y} \in H$ with $\bar{y} = g(\bar{x})$, by )2.3) and (2.4), the covering constant for $g$ at point $(\bar{x}, \bar{y})$

$$\hat{\alpha}(g, \bar{x}, \bar{y}) = \sup_{\eta > 0} \inf\{\|z\| : z \in \widehat{D}^* g(x)(y), x \in \mathbb{B}_X(\bar{x}, \eta), g(x) \in \mathbb{B}_Y(\bar{y}, \eta), \|y\| = 1\}. \quad (2.6)$$

In this paper, we concentrate to Euclidean spaces. We first review the results about Fréchet derivative and Mordukhovich derivative obtained in [15].

Let $n \geq 1$ and let $(\mathbb{R}^n, \|\cdot\|)$ be the standard $n$-d Euclidean space with the ordinal Hilbert $L_2$-norm and row vectors. Let $\theta$ denote the origin of $\mathbb{R}^n$. Let $m, n \geq 1$ and let $f: \mathbb{R}^n \to \mathbb{R}^m$ be a single-valued mapping with the following representation.

$$f((x_1, x_2, \ldots, x_n)) = \big(f_1(x_1, x_2, \ldots, x_n), f_2(x_1, x_2, \ldots, x_n), \ldots, f_m(x_1, x_2, \ldots, x_n)\big), \text{ for } (x_1, x_2, \ldots, x_n) \in \mathbb{R}^n.$$

Where, for $i = 1, 2, \ldots, m$, $f_i(x_1, x_2, \ldots, x_n): \mathbb{R}^n \to \mathbb{R}$ is a real valued multivariable function defined on $\mathbb{R}^n$ with variables $x_1, x_2, \ldots, x_n$. We write the above equation as $f = (f_1, f_2, \ldots, f_m)$. Let $(z_1, z_2, \ldots, z_n) \in \mathbb{R}^n$. For $i = 1, 2, \ldots, n$, the partial derivative of $f_i$ with respect to the variable $x_j$ at $(z_1, z_2, \ldots, z_n)$ is

$$\frac{\partial f_i}{\partial x_j}(z_1, z_2, \ldots, z_n), \text{ for every } j = 1, 2, \ldots, n.$$

Then, for any $i = 1, 2, \ldots, m$, the existence of $f_i(z_1, z_2, \ldots, z_n)$ means that,

$$\lim_{x_j \to z_j} \frac{f_i(x_1, x_2, \ldots, x_n) - f_i(z_1, z_2, \ldots, z_n) - \frac{\partial f_i}{\partial x_j}(z_1, z_2, \ldots, z_n)(x_j - z_j)}{x_j - z_j} = 0, \text{ for } j = 1, 2, \ldots, n.$$

This limit is equivalent to, for any $i = 1, 2, \ldots, m$,

$$\lim_{x_j \to z_j} \frac{f_i(x_1, x_2, \ldots, x_n) - f_i(z_1, z_2, \ldots, z_n)}{x_j - z_j} = -\frac{\partial f_i}{\partial x_j}(z_1, z_2, \ldots, z_n), \text{ for } j = 1, 2, \ldots, n.$$

Furthermore, if $f$ satisfies the following conditions

$$\lim_{x \to z} \frac{f_i(x_1, x_2, \ldots, x_n) - f_i(z_1, z_2, \ldots, z_n) - \left(\sum_{j=1}^{n} \frac{\partial f_i}{\partial x_j}(z_1, z_2, \ldots, z_n)(x_j - z_j)\right)}{\|x - z\|} = 0, \text{ for } i = 1, 2, \ldots, m,$$

then, $f$ is differentiable and has the linear approximation at point $z = (z_1, z_2, \ldots, z_n)$, which is the first order of the Taylor polynomial of the mapping $f$ from $\mathbb{R}^n$ to $\mathbb{R}^m$. We have some sufficient conditions for $f$ to have the linear approximation at point $z$.

For every $i = 1, 2, \ldots, m$, the second order partial derivative of $f_i$ at point $z = (z_1, z_2, \ldots, z_n)$ with respect to $x_j$ and $x_k$ is

$$\frac{\partial^2 f_i}{\partial x_j \partial x_k}(z_1, z_2, \ldots, z_n), \text{ for } j, k = 1, 2, \ldots, n.$$

**Fact 3.1 in [15]**. *Let $f = (f_1, f_2, \ldots, f_m)$. $: \mathbb{R}^n \to \mathbb{R}^m$ be a single-valued mapping. Let $z \in \mathbb{R}^n$. If there is a ball $B$ with radius $r > 0$ and centered at $z$ such that, for every $i = 1, 2, \ldots, m$, the real valued function $f_i: \mathbb{R}^n \to \mathbb{R}$ is twice differentiable (all second partial derivatives of $f_i$ exist) at every point $y \in B$, that is*

$$\frac{\partial^2 f_i}{\partial x_j \partial x_k}(z_1, z_2, \ldots, z_n) \text{ exists for any } z \in B, \text{ for } j, k = 1, 2, \ldots, n.$$

*then $f$ has the linear approximation at this point $z$.*

**Theorem 3.2 in [15]**. *Let $f = (f_1, f_2, \ldots, f_m): \mathbb{R}^n \to \mathbb{R}^m$ be a single-valued mapping. Let $z = (z_1, z_2, \ldots, z_n) \in \mathbb{R}^n$. Suppose that, for every $i = 1, 2, \ldots, m$, $\frac{\partial f_i}{\partial x_j}(z_1, z_2, \ldots, z_n)$ exists, for every $j = 1, 2, \ldots, n$ and $f$ has the linear approximation at point $z$. Then,*

(a) *$f$ is Fréchet differentiable at $z$ and the Fréchet derivative of $f$ at $z$ is the following $n \times m$ matrix,*

$$\nabla f(z) = \begin{pmatrix} \frac{\partial f_1}{\partial x_1}(z_1, z_2, \ldots, z_n) & \cdots & \frac{\partial f_m}{\partial x_1}(z_1, z_2, \ldots, z_n) \\ \vdots & \ddots & \vdots \\ \frac{\partial f_1}{\partial x_n}(z_1, z_2, \ldots, z_n) & \cdots & \frac{\partial f_m}{\partial x_n}(z_1, z_2, \ldots, z_n) \end{pmatrix},$$

(b) *$f$ is Mordukhovich differentiable at point $z$ that is the <span style="color:red">Jacobian matrix of $f$ at $z$</span>*

$$\widehat{D}^* f(x) = \nabla f(z)^T = \begin{pmatrix} \frac{\partial f_1}{\partial x_1}(z_1, z_2, \ldots, z_n) & \cdots & \frac{\partial f_1}{\partial x_n}(z_1, z_2, \ldots, z_n) \\ \vdots & \ddots & \vdots \\ \frac{\partial f_m}{\partial x_1}(z_1, z_2, \ldots, z_n) & \cdots & \frac{\partial f_m}{\partial x_n}(z_1, z_2, \ldots, z_n) \end{pmatrix}.$$

**Proposition 3.3 in [15]**. *Let $f = (f_1, f_2, \ldots, f_m): \mathbb{R}^n \to \mathbb{R}^m$ be a single-valued mapping.*

(a) *Suppose that, for every $i = 1, 2, \ldots, m$, $f_i(x_1, x_2, \ldots, x_n)$ is a polynomial function with respect to $x_1, x_2, \ldots, x_n$, then $f$ is Fréchet differentiable at every point in $\mathbb{R}^n$, and $f$ is Mordukhovich differentiable on $\mathbb{R}^n$;*

(b) *Suppose that, for every $i = 1, 2, \ldots, m$, $f_i(x_1, x_2, \ldots, x_n)$ is a rational function with respect to $x_1, x_2, \ldots, x_n$. Let $z \in \mathbb{R}^n$. If $z$ is not a zero point for the denominator of every $f_i$, then $f$ is*

*Fréchet differentiable at z; and therefore, f is Mordukhovich differentiable at z.*

## 3. A Norm preserving Mapping in $\mathbb{R}^2$

In this section, we consider the single-valued mapping $f: \mathbb{R}^2 \to \mathbb{R}^2$ defined by (1.1), which is studied in Examples 2 in [4] and 4.2 in [2]. We will calculate the Fréchet and Mordukhovich derivatives of $f$, by which we will find the exact covering constant for $f$. Recall that $f$ is defined by

$$f((x_1, x_2)) = \left( \frac{x_1^2 - x_2^2}{\sqrt{x_1^2 + x_2^2}}, \frac{2x_1 x_2}{\sqrt{x_1^2 + x_2^2}} \right), \text{ for } (x_1, x_2) \in \mathbb{R}^2 \setminus \{\theta\} \text{ with } f(\theta) = \theta.$$

We have that $\|f(\theta)\| = \|\theta\|$. For any $x = (x_1, x_2) \in \mathbb{R}^2$ with $(x_1, x_2) \neq \theta$, we calculate

$$\|f((x_1, x_2))\|^2 = \left( \frac{x_1^2 - x_2^2}{\sqrt{x_1^2 + x_2^2}} \right)^2 + \left( \frac{2x_1 x_2}{\sqrt{x_1^2 + x_2^2}} \right)^2 = x_1^2 + x_2^2 = \|(x_1, x_2)\|^2.$$

This implies that the mapping $f: \mathbb{R}^2 \to \mathbb{R}^2$ is norm preserving. That is,

$$\|f(x)\| = \|x\|, \text{ for any } x \in \mathbb{R}^2.$$

By the norm preserving property of $f$, it follows immediately that $f$ is a continuous mapping on $\mathbb{R}^2$, which means that $f$ is continuous at origin $\theta$ in $\mathbb{R}^2$. Next, we prove that at the origin $\theta$ in $\mathbb{R}^2$, $f$ is neither Fréchet differentiable, nor Mordukhovich differentiable.

**Lemma 3.1.** *$f$ is not Fréchet differentiable at $\theta$.*

*Proof.* Assume, by the way of contradiction, that $f$ is Fréchet differentiable at $\theta$ and $\nabla f(\theta): \mathbb{R}^2 \to \mathbb{R}^2$ exists, which is a linear and continuous mapping on $\mathbb{R}^2$. Then, $\nabla f(\theta)$ must be precisely represented by a real $2 \times 2$ matrix $\begin{pmatrix} a & b \\ c & d \end{pmatrix}$, for some real numbers $a$, $b$, $c$, and $d$.

Assume $a \neq 0$. In this case, we take a special direction in the following limit for $u \to \theta$ by $u_1 = -at$ with $t \downarrow 0$, $u_2 = 0$. By the assumption that $\nabla f(\theta) = \begin{pmatrix} a & b \\ c & d \end{pmatrix}$ and by $f(\theta) = \theta$, it follows that

$$\lim_{u=(u_1, u_2) \to \theta} \frac{f(u) - f(\theta) - \nabla f(z)(u - \theta)}{\|u - z\|}$$

$$= \lim_{t \downarrow 0} \frac{\left( \frac{a^2 t^2}{\sqrt{a^2 t^2}}, 0 \right) - (-at, 0)\begin{pmatrix} a & b \\ c & d \end{pmatrix}}{\|(-at, 0)\|}$$

$$= \lim_{t \downarrow 0} \frac{(|a|t, 0) + (a^2 t, abt)}{|a|t}$$

$$= (1 + |a|, (\text{sign}a)b)$$

$$\neq \theta.$$

This contradicts to the assumption that $\nabla f(\theta) = \begin{pmatrix} a & b \\ c & d \end{pmatrix}$.

Next, suppose $a = 0$. In this case, we take a special direction in the limit (5.3) for $u \to \theta$ by $u_1 \downarrow 0$, $u_2 = 0$, $u_3 = z_3$ and $u_4 = z_4$. It follows that

$$\lim_{u=(u_1,u_2) \to \theta} \frac{f(u)-f(\theta)-\nabla f(z)(u-\theta)}{\|u-z\|}$$

$$= \lim_{u_1 \downarrow 0} \frac{\left(\frac{u_1^2}{\sqrt{u_1^2}}, 0\right) - (u_1, 0)\begin{pmatrix} 0 & b \\ c & d \end{pmatrix}}{\|(u_1, 0)\|}$$

$$= \lim_{u_1 \downarrow 0} \frac{(1, 0) - (0, bu_1)}{u_1}$$

$$= (1, b)$$

$$\neq \theta.$$

This contradicts the assumption that $\nabla f(\theta) = \begin{pmatrix} a & b \\ c & d \end{pmatrix}$. It completes the proof. $\square$

**Proposition 3.2.** $f$ is not Mordukhovich differentiable at $\theta$. More precisely speaking, we have

$$\widehat{D}^* f(\theta)(y) = \emptyset, \quad \text{for any } y \in \mathbb{R}^2 \setminus \{\theta\}.$$

*Proof.* Let $y = (y_1, y_2) \in \mathbb{R}^2$ and $x = (x_1, x_2) \in \mathbb{R}^2$, suppose that $x \in \widehat{D}^* f(\theta)(y)$. We calculate the following limits.

$$0 \geq \limsup_{u \to \theta} \frac{\langle (x_1, x_2), (u_1, u_2) - \theta \rangle - \langle (y_1, y_2), f(u) - f(\theta) \rangle}{\|u - \theta\| + \|f(u) - f(\theta)\|}$$

$$= \limsup_{u \to \theta} \frac{\langle (x_1, x_2), (u_1, u_2) \rangle - \langle (y_1, y_2), \left( \frac{u_1^2 - u_2^2}{\sqrt{u_1^2 + u_2^2}}, \frac{2u_1 u_2}{\sqrt{u_1^2 + u_2^2}} \right) \rangle}{\|u\| + \|f(u)\|}$$

$$= \limsup_{u \to \theta} \frac{x_1 u_1 + x_2 u_2 - y_1 \frac{u_1^2 - u_2^2}{\sqrt{u_1^2 + u_2^2}} - y_2 \frac{2u_1 u_2}{\sqrt{u_1^2 + u_2^2}}}{2\|u\|}. \tag{3.1}$$

Case 1. $x_1 > y_1$. In this case, we take a special direction in the limit (3.1) for $u \to \theta$ by $u_2 = 0$ and $u_1 \downarrow 0$. It follows that

$$\limsup_{u \to \theta} \frac{x_1 u_1 + x_2 u_2 - y_1 \frac{u_1^2 - u_2^2}{\sqrt{u_1^2 + u_2^2}} - y_2 \frac{2u_1 u_2}{\sqrt{u_1^2 + u_2^2}}}{2\|u\|}$$

$$\geq \lim_{u_2=0 \text{ and } u_1 \downarrow 0} \frac{x_1 u_1 - y_1 u_1}{2|u_1|}$$

$$= \lim_{u_2=0 \text{ and } u_1 \downarrow 0} \left( \frac{x_1}{2} - \frac{y_1}{2} \right)$$

$$= \frac{x_1}{2} - \frac{y_1}{2} > 0.$$

This implies that

$$x_1 > y_1 \quad \Longrightarrow \quad x \notin \widehat{D}^* f(\theta)(y).$$

Case 2. $-x_1 > y_1$. In this case, we take a special direction in the limit (3.1) for $u \to \theta$ by $u_2 = 0$ and $u_1 \uparrow 0$. We have

$$\limsup_{u \to \theta} \frac{x_1 u_1 + x_2 u_2 - y_1 \frac{u_1^2 - u_2^2}{\sqrt{u_1^2 + u_2^2}} - y_2 \frac{2 u_1 u_2}{\sqrt{u_1^2 + u_2^2}}}{2\|u\|}$$

$$\geq \lim_{u_2 = 0 \text{ and } u_1 \uparrow 0} \frac{x_1 u_1 - y_1 |u_1|}{2|u_1|}$$

$$= \lim_{u_2 = 0 \text{ and } u_1 \uparrow 0} \left( \frac{-x_1}{2} - \frac{y_1}{2} \right)$$

$$= -\frac{x_1}{2} - \frac{y_1}{2} > 0.$$

This implies that

$$-x_1 > y_1 \quad \Longrightarrow \quad x \notin \widehat{D}^* f(\theta)(y).$$

Case 3. $x_1 + x_2 > \sqrt{2} y_2$. In this case, we take a special direction in the limit (3.1) for $u \to \theta$ by $u_1 = u_2$ and $u_1 \downarrow 0$. We have

$$\limsup_{u \to \theta} \frac{x_1 u_1 + x_2 u_2 - y_1 \frac{u_1^2 - u_2^2}{\sqrt{u_1^2 + u_2^2}} - y_2 \frac{2 u_1 u_2}{\sqrt{u_1^2 + u_2^2}}}{2\|u\|}$$

$$\geq \lim_{u_1 = u_2 \text{ and } u_1 \downarrow 0} \frac{x_1 u_1 + x_2 u_1 - \frac{2 y_2 u_1}{\sqrt{2}}}{2\sqrt{2} u_1}$$

$$= \lim_{u_1 = u_2 \text{ and } u_1 \downarrow 0} \frac{x_1 + x_2 - \sqrt{2} y_2}{2\sqrt{2}}$$

$$= \frac{x_1 + x_2 - \sqrt{2} y_2}{2\sqrt{2}} > 0.$$

This implies that

$$x_1 + x_2 > \sqrt{2} y_2 \quad \Longrightarrow \quad x \notin \widehat{D}^* f(\theta)(y).$$

Case 4. $-x_1 - x_2 > \sqrt{2} y_2$. In this case, we take a special direction in the limit (3.1) for $u \to \theta$ by $u_1 = u_2$ and $u_1 \uparrow 0$. We have

$$\limsup_{u \to \theta} \frac{x_1 u_1 + x_2 u_2 - y_1 \frac{u_1^2 - u_2^2}{\sqrt{u_1^2 + u_2^2}} - y_2 \frac{2 u_1 u_2}{\sqrt{u_1^2 + u_2^2}}}{2\|u\|}$$

$$\geq \lim_{u_1 = u_2 \text{ and } u_1 \uparrow 0} \frac{x_1 u_1 + x_2 u_2 - \frac{2 y_2 |u_1|}{\sqrt{2}}}{2\sqrt{2} |u_1|}$$

$$= \lim_{u_1=u_2 \text{ and } u_1\uparrow 0} \frac{-x_1-x_2-\sqrt{2}y_2}{2\sqrt{2}}$$

$$= \frac{-x_1-x_2-\sqrt{2}y_2}{2\sqrt{2}} > 0.$$

This implies that

$$-x_1 - x_2 > \sqrt{2}y_2 \quad \Rightarrow \quad x \notin \widehat{D}^*f(\theta)(y).$$

Case 5. $-x_1 + x_2 > -\sqrt{2}y_2$. In this case, we take a special direction in the limit (3.1) for $u \to \theta$ by $u_1 = -u_2$ and $u_2 \downarrow 0$. It follows that

$$\limsup_{u \to \theta} \frac{x_1u_1 + x_2u_2 - y_1 \frac{u_1^2 - u_2^2}{\sqrt{u_1^2+u_2^2}} - y_2 \frac{2u_1u_2}{\sqrt{u_1^2+u_2^2}}}{2\|u\|}$$

$$\geq \lim_{u_1=-u_2 \text{ and } u_2\downarrow 0} \frac{x_1u_1 + x_2u_2 + \sqrt{2}y_2u_2}{2\sqrt{2}|u_2|}$$

$$= \lim_{u_1=-u_2 \text{ and } u_2\downarrow 0} \frac{-x_1u_2 + x_2u_2 + \sqrt{2}y_2u_2}{2\sqrt{2}u_2}$$

$$= \lim_{u_1=-u_2 \text{ and } u_2\downarrow 0} \left( \frac{-x_1}{2\sqrt{2}} + \frac{x_2}{2\sqrt{2}} + \frac{\sqrt{2}y_2}{2\sqrt{2}} \right)$$

$$= \frac{-x_1}{2\sqrt{2}} + \frac{x_2}{2\sqrt{2}} + \frac{\sqrt{2}y_2}{2\sqrt{2}} > 0.$$

This implies that

$$-x_1 + x_2 > -\sqrt{2}y_2 \quad \Rightarrow \quad x \notin \widehat{D}^*f(\theta)(y).$$

Case 6. $x_1 - x_2 > -\sqrt{2}y_2$. In this case, we take a special direction in the limit (3.1) for $u \to \theta$ by $u_1 = -u_2$ and $u_2 \uparrow 0$. It follows that

$$\limsup_{u \to \theta} \frac{x_1u_1 + x_2u_2 - y_1 \frac{u_1^2 - u_2^2}{\sqrt{u_1^2+u_2^2}} - y_2 \frac{2u_1u_2}{\sqrt{u_1^2+u_2^2}}}{2\|u\|}$$

$$\geq \lim_{u_1=-u_2 \text{ and } u_2\uparrow 0} \frac{x_1u_1 + x_2u_2 + \sqrt{2}y_2|u_2|}{2\sqrt{2}|u_2|}$$

$$= \lim_{u_1=-u_2 \text{ and } u_2\uparrow 0} \frac{-x_1u_2 + x_2u_2 + \sqrt{2}y_2|u_2|}{2\sqrt{2}|u_2|}$$

$$= \lim_{u_1=-u_2 \text{ and } u_2\uparrow 0} \left( \frac{x_1}{2\sqrt{2}} - \frac{x_2}{2\sqrt{2}} + \frac{\sqrt{2}y_2}{2\sqrt{2}} \right)$$

$$= \frac{x_1}{2\sqrt{2}} - \frac{x_2}{2\sqrt{2}} + \frac{\sqrt{2}y_2}{2\sqrt{2}} > 0.$$

This implies that

$$x_1 - x_2 > -\sqrt{2}y_2 \quad \Rightarrow \quad x \notin \widehat{D}^*f(\theta)(y).$$

By summarizing the above 6 cases, we have that the following inequalities are necessary conditions for $x \in \widehat{D}^*f(\theta)(y)$.

(i) $x_1 \leq y_1$

(ii) $-x_1 \leq y_1$

(iii) $x_1 + x_2 \leq \sqrt{2}y_2$

(iv) $-x_1 - x_2 \leq \sqrt{2}y_2$

(v) $-x_1 + x_2 \leq -\sqrt{2}y_2$

(vi) $x_1 - x_2 \leq -\sqrt{2}y_2$

From the above inequalities (i–vi), we obtain more precise necessary conditions for $x \in \widehat{D}^*f(\theta)(y)$.

By (iii–vi), we have

(I) $\qquad\qquad\qquad x \in \widehat{D}^*f(\theta)(y) \quad \Rightarrow \quad y_2 = 0.$

Substituting $y_2 = 0$ into (v–vi), we get

(II) $\qquad\qquad\qquad x \in \widehat{D}^*f(\theta)(y) \quad \Rightarrow \quad x_1 - x_2 = 0$, that is $x_1 = x_2$.

Substituting $y_2 = 0$ and $x_1 = x_2$ into (iii–iv), we get

(III) $\qquad\qquad\qquad x \in \widehat{D}^*f(\theta)(y) \quad \Rightarrow \quad x_1 = x_2 = 0.$

Substituting $x_1 = 0$ into (i) or (ii), we get

(IV) $\qquad\qquad\qquad x \in \widehat{D}^*f(\theta)(y) \quad \Rightarrow \quad y_1 \geq 0.$

Under the necessary conditions (I–III), assume $y_1 > 0$. In this case, we take a special direction in the limit (3.1) for $u \to \theta$ by $u_1 = 0$ and $u_2 \uparrow 0$. We have

$$\limsup_{u \to \theta} \frac{x_1 u_1 + x_2 u_2 - y_1 \frac{u_1^2 - u_2^2}{\sqrt{u_1^2 + u_2^2}} - y_2 \frac{2u_1 u_2}{\sqrt{u_1^2 + u_2^2}}}{2\|u\|}$$

$$= \limsup_{u \to \theta} \frac{-y_1 \frac{u_1^2 - u_2^2}{\sqrt{u_1^2 + u_2^2}}}{2\|u\|}$$

$$\geq \lim_{u_1 = 0 \text{ and } u_2 \uparrow 0} \frac{y_1 |u_2|}{2|u_2|}$$

$$= \frac{y_1}{2} > 0.$$

This implies that

$$y_1 > 0 \quad \Longrightarrow \quad x \notin \widehat{D}^* f(\theta)(y).$$

By the necessary condition (IV), this implies that

(V) $$\qquad\qquad x \in \widehat{D}^* f(\theta)(y) \quad \Longrightarrow \quad y_1 = 0.$$

Under the necessary conditions (I–V), we obtain that, for $y = (y_1, y_2) \in \mathbb{R}^2$ and $x = (x_1, x_2) \in \mathbb{R}^2$,

$$x \in \widehat{D}^* f(\theta)(y) \quad \Longrightarrow \quad x = y = \theta.$$

This proposition is proved. □

Next, we show that $f$ is Fréchet differentiable and Mordukhovich differentiable on $\mathbb{R}^2 \setminus \{\theta\}$.

**Theorem 3.3.** *Let* $z = (z_1, z_2) \in \mathbb{R}^2 \setminus \{\theta\}$. *Then*

(a) *$f$ is Fréchet differentiable at $z$ and*

$$\nabla f(z) = \begin{pmatrix} \dfrac{(z_1^2 + 3z_2^2)z_1}{(z_1^2+z_2^2)\sqrt{z_1^2+z_2^2}} & \dfrac{2z_2^2 z_2}{(z_1^2+z_2^2)\sqrt{z_1^2+z_2^2}} \\ \dfrac{-(3z_1^2+z_2^2)z_2}{(z_1^2+z_2^2)\sqrt{z_1^2+z_2^2}} & \dfrac{2z_1^2 z_1}{(z_1^2+z_2^2)\sqrt{z_1^2+z_2^2}} \end{pmatrix}.$$

(b) *$f$ is Mordukhovich differentiable at $z$ and*

$$\widehat{D}^* f(z) = \begin{pmatrix} \dfrac{(z_1^2+3z_2^2)z_1}{(z_1^2+z_2^2)\sqrt{z_1^2+z_2^2}} & \dfrac{-(3z_1^2+z_2^2)z_2}{(z_1^2+z_2^2)\sqrt{z_1^2+z_2^2}} \\ \dfrac{2z_2^2 z_2}{(z_1^2+z_2^2)\sqrt{z_1^2+z_2^2}} & \dfrac{2z_1^2 z_1}{(z_1^2+z_2^2)\sqrt{z_1^2+z_2^2}} \end{pmatrix}.$$

*However, more precisely speaking, for* $x = (x_1, x_2), y = (y_1, y_2) \in \mathbb{R}^2$, *if* $x = \widehat{D}^* f(z)(y)$, *then*

$$x_1 = y_1 \frac{(z_1^2+3z_2^2)z_1}{(z_1^2+z_2^2)\sqrt{z_1^2+z_2^2}} + y_2 \frac{2z_2^2 z_2}{(z_1^2+z_2^2)\sqrt{z_1^2+z_2^2}},$$

$$x_2 = y_1 \frac{(-3z_1^2-z_2^2)z_2}{(z_1^2+z_2^2)\sqrt{z_1^2+z_2^2}} + y_2 \frac{2z_1^2 z_1}{(z_1^2+z_2^2)\sqrt{z_1^2+z_2^2}}.$$

*Proof.* Proof of (a). Let $f = (f_1, f_2)$ with

$$f_1(x_1, x_2) = \frac{x_1^2 - x_2^2}{\sqrt{x_1^2+x_2^2}}, \text{ for } (x_1, x_2) \in \mathbb{R}^2 \setminus \{\theta\} \text{ with } f_1(\theta) = 0,$$

and $$f_2(x_1, x_2) = \frac{2x_1 x_2}{\sqrt{x_1^2+x_2^2}}, \text{ for } (x_1, x_2) \in \mathbb{R}^2 \setminus \{\theta\} \text{ with } f_2(\theta) = 0.$$

Let $z = (z_1, z_2) \in \mathbb{R}^2 \setminus \{\theta\}$. There is $r > 0$ such that for any $x = (x_1, x_2) \in \mathbb{R}^2$,

$$\|x - z\| < r \implies x \neq \theta.$$

This implies that both $f_1$ and $f_2$ are twice differentiable on the ball with radius $r > 0$ and centered at $z$. Then by Lemma 3.1 and Theorem 3.2 in [15], $f$ is Fréchet differentiable at $z$ and the Fréchet derivative of $f$ at $z$ is the following 2×2 matrix,

$$\nabla f(z) = \begin{pmatrix} \frac{(z_1^2+3z_2^2)z_1}{(z_1^2+z_2^2)\sqrt{z_1^2+z_2^2}} & \frac{2z_2^2 z_2}{(z_1^2+z_2^2)\sqrt{z_1^2+z_2^2}} \\ \frac{-(3z_1^2+z_2^2)z_2}{(z_1^2+z_2^2)\sqrt{z_1^2+z_2^2}} & \frac{2z_1^2 z_1}{(z_1^2+z_2^2)\sqrt{z_1^2+z_2^2}} \end{pmatrix}.$$

Proof of part (b). By Theorem 1.38 in [15] and Part (a), we have that $f$ is Mordukhovich differentiable at z and it satisfies that

$$\widehat{D}^* f(z)(y) = (\nabla f(x))^*(y), \text{ for any } y \in \mathbb{R}^2.$$

By the representation of $\nabla f(x)$ in Part (a), this gives the solution of $\widehat{D}^* f(z)$ listed in Part (b). □

**Proposition 3.4.** *Let* $z = (z_1, z_2) \in \mathbb{R}^2 \setminus \{\theta\}$ *and* $x = (x_1, x_2), y = (y_1, y_2) \in \mathbb{R}^2$, *if* $x = \widehat{D}^* f(z)(y)$, *then x and y satisfy the following conditions.*

(i) $\qquad x_1^2 + x_2^2 = y_1^2 + y_2^2 + 12\left(y_1 \frac{z_1 z_2}{z_1^2+z_2^2} - y_2 \frac{z_1^2-z_2^2}{2(z_1^2+z_2^2)}\right)^2.$

(ii) $\qquad \|x\| \geq \|y\|.$

(iii) $\qquad y_1 z_1 z_2 = y_2 \frac{z_1^2-z_2^2}{2} \implies \|x\| = \|y\|.$

*Proof.* We only prove (i). By Proposition 3.6 and $x \in \widehat{D}^* f(z)(y)$, we have

$$x_1 = y_1 \left(\frac{(z_1^2+3z_2^2)z_1}{(z_1^2+z_2^2)\sqrt{z_1^2+z_2^2}}\right) + y_2 \left(\frac{2z_2^2 z_2}{(z_1^2+z_2^2)\sqrt{z_1^2+z_2^2}}\right),$$

and $\qquad x_2 = y_1 \left(\frac{-(3z_1^2+z_2^2)z_2}{(z_1^2+z_2^2)\sqrt{z_1^2+z_2^2}}\right) + y_2 \left(\frac{2z_1^2 z_1}{(z_1^2+z_2^2)\sqrt{z_1^2+z_2^2}}\right).$

Then, we calculate

$$x_1^2 + x_2^2 = \left(y_1 \left(\frac{(z_1^2+3z_2^2)z_1}{(z_1^2+z_2^2)\sqrt{z_1^2+z_2^2}}\right) + y_2 \left(\frac{2z_2^2 z_2}{(z_1^2+z_2^2)\sqrt{z_1^2+z_2^2}}\right)\right)^2 + \left(y_1 \left(\frac{-(3z_1^2+z_2^2)z_2}{(z_1^2+z_2^2)\sqrt{z_1^2+z_2^2}}\right) + y_2 \left(\frac{2z_1^2 z_1}{(z_1^2+z_2^2)\sqrt{z_1^2+z_2^2}}\right)\right)^2$$

$$= y_1^2\left(\left(\frac{(z_1^2+3z_2^2)z_1}{(z_1^2+z_2^2)\sqrt{z_1^2+z_2^2}}\right)^2 + \left(\frac{-(3z_1^2+z_2^2)z_2}{(z_1^2+z_2^2)\sqrt{z_1^2+z_2^2}}\right)^2\right) + y_2^2\left(\left(\frac{2z_2^2z_2}{(z_1^2+z_2^2)\sqrt{z_1^2+z_2^2}}\right)^2 + \left(\frac{2z_1^2z_1}{(z_1^2+z_2^2)\sqrt{z_1^2+z_2^2}}\right)^2\right)$$

$$+ 2y_1y_2\left(\left(\frac{-(3z_1^2+z_2^2)z_2}{(z_1^2+z_2^2)\sqrt{z_1^2+z_2^2}}\right)\left(\frac{2z_2^2z_2}{(z_1^2+z_2^2)\sqrt{z_1^2+z_2^2}}\right) + \left(\frac{-(3z_1^2+z_2^2)z_2}{(z_1^2+z_2^2)\sqrt{z_1^2+z_2^2}}\right)\left(\frac{2z_1^2z_1}{(z_1^2+z_2^2)\sqrt{z_1^2+z_2^2}}\right)\right)$$

$$= y_1^2\left(\frac{(z_1^2-z_2^2)^2}{(z_1^2+z_2^2)^2} - \frac{4(z_1^2-z_2^2)^2}{(z_1^2+z_2^2)^2} + \frac{4(z_1^2+z_2^2)}{z_1^2+z_2^2}\right) + 2y_1y_2\left(\frac{2(z_1^2-z_2^2)z_1z_2}{(z_1^2+z_2^2)^2} - \frac{4(z_1^2-z_2^2)z_1z_2}{(z_1^2+z_2^2)^2} - \frac{4(z_1^2-z_2^2)z_1z_2}{(z_1^2+z_2^2)^2}\right)$$

$$+ y_2^2\left(4\left(\frac{z_1^4+2z_1^2z_2^2+z_2^4}{(z_1^2+z_2^2)^2} - \frac{3z_1^2z_2^2}{(z_1^2+z_2^2)^2}\right)\right)$$

$$= y_1^2\left(\frac{z_1^4+14z_1^2z_2^2+z_2^4}{(z_1^2+z_2^2)^2}\right) + 2y_1y_2\left(-\frac{6(z_1^2-z_2^2)z_1z_2}{(z_1^2+z_2^2)^2}\right) + y_2^2\left(4\left(\frac{z_1^4+2z_1^2z_2^2+z_2^4}{(z_1^2+z_2^2)^2} - \frac{3z_1^2z_2^2}{(z_1^2+z_2^2)^2}\right)\right)$$

$$= y_1^2\left(1 + \frac{12z_1^2z_2^2}{(z_1^2+z_2^2)^2}\right) + 2y_1y_2\left(-\frac{6(z_1^2-z_2^2)z_1z_2}{(z_1^2+z_2^2)^2}\right) + y_2^2\left(4\left(\frac{(z_1^2-z_2^2)^2}{(z_1^2+z_2^2)^2} + \frac{z_1^2z_2^2}{(z_1^2+z_2^2)^2}\right)\right)$$

$$= y_1^2 + y_2^2 + y_1^2\frac{12z_1^2z_2^2}{(z_1^2+z_2^2)^2} - 4y_1y_2\left(\frac{3(z_1^2-z_2^2)z_1z_2}{(z_1^2+z_2^2)^2}\right) + 4y_2^2\left(\frac{(z_1^2-z_2^2)^2}{(z_1^2+z_2^2)^2} + \frac{z_1^2z_2^2}{(z_1^2+z_2^2)^2} - \frac{1}{4}\right)$$

$$= y_1^2 + y_2^2 + y_1^2\frac{12z_1^2z_2^2}{(z_1^2+z_2^2)^2} - 4y_1y_2\left(\frac{3(z_1^2-z_2^2)z_1z_2}{(z_1^2+z_2^2)^2}\right) + 4y_2^2\left(\frac{4(z_1^2-z_2^2)^2+4z_1^2z_2^2-(z_1^2+z_2^2)^2}{4(z_1^2+z_2^2)^2}\right)$$

$$= y_1^2 + y_2^2 + y_1^2\frac{12z_1^2z_2^2}{(z_1^2+z_2^2)^2} - 4y_1y_2\left(\frac{3(z_1^2-z_2^2)z_1z_2}{(z_1^2+z_2^2)^2}\right) + 4y_2^2\left(\frac{3(z_1^4+z_2^4)-6z_1^2z_2^2}{4(z_1^2+z_2^2)^2}\right)$$

$$= y_1^2 + y_2^2 + 12\left(y_1^2\frac{z_1^2z_2^2}{(z_1^2+z_2^2)^2} - y_1y_2\frac{(z_1^2-z_2^2)z_1z_2}{(z_1^2+z_2^2)^2} + y_2^2\frac{(z_1^2-z_2^2)^2}{4(z_1^2+z_2^2)^2}\right)$$

$$= y_1^2 + y_2^2 + 12\left(y_1\frac{z_1z_2}{z_1^2+z_2^2} - y_2\frac{z_1^2-z_2^2}{2(z_1^2+z_2^2)}\right)^2. \qquad \square$$

**Theorem 3.5.** *Let $f: \mathbb{R}^2 \to \mathbb{R}^2$ be defined by (4.1). Then, we have*

$$\hat{\alpha}(f, \bar{z}, f(\bar{z})) = 1, \text{ for any } \bar{z} \in \mathbb{R}^2.$$

*Proof.* For any $\bar{z}, \bar{w} \in \mathbb{R}^2$ with $\bar{w} = f(\bar{z})$ and $\eta > 0$, let $\mathbb{B}(\bar{z}, \eta)$ and $\mathbb{B}(\bar{w}, \eta)$ denote the closed ball in $\mathbb{R}^2$ with radius $\eta$ and centered at point $\bar{z}$ and $\bar{w}$, respectively. By (2.4), we have

$$\hat{\alpha}(f, \bar{z}, \bar{w}) = \sup_{\eta>0} \inf\{\|x\|: x \in \widehat{D}^*f(z)(y), z \in \mathbb{B}(\bar{z}, \eta), f(z) \in \mathbb{B}(\bar{w}, \eta), \|y\| = 1\}$$

$$= \sup_{\eta>0} \inf\{\|x\|: x \in \widehat{D}^*f(z)(y), z \in \mathbb{B}(\bar{z}, \eta), f(z) \in \mathbb{B}(\bar{w}, \eta), \|y\| = 1\}.$$

In case $\theta \in \mathbb{B}(\bar{z}, \eta)$ and $f(\theta) = \theta \in \mathbb{B}(\bar{w}, \eta)$, by Proposition 4.4, we have that

$$\widehat{D}^* f(\theta)(y) = \emptyset, \text{ for any } y \in \mathbb{R}^2 \text{ with } \|y\| = 1.$$

This implies that in the above limit set in, there are some point $z \in \mathbb{B}(\bar{z}, \eta)$ and $y \in \mathbb{R}^2$ with $\|y\| = 1$ such that $\widehat{D}^* f(z)(y) = \emptyset$. Since this Hilbert space $\mathbb{R}^2$ is Asplund and $f: \mathbb{R}^2 \to \mathbb{R}^2$ defined by (1.1) is continuous on $\mathbb{R}^2$, it yields that, for every $\eta > 0$, $f$ is Fréchet differentiable on a dense subset of $\mathbb{B}(\bar{z}, \eta)$. Then, by Theorem 1.38 in [15], $f$ is Mordukhovich differentiable on a dense subset of $\mathbb{B}(\bar{z}, \eta)$, for every $\eta > 0$, which implies that

$$\{\|x\|: x \in \widehat{D}^* f(z)(y), z \in \mathbb{B}(\bar{z}, \eta), f(z) \in \mathbb{B}(\bar{w}, \eta), \|y\| = 1\} \neq \emptyset, \text{ for every } \eta > 0.$$

However, if $\widehat{D}^* f(z)(y) \neq \emptyset$, for $z \in \mathbb{B}(\bar{z}, \eta)$ and $y \in \mathbb{R}^2$ with $\|y\| = 1$, then, for any $x \in \widehat{D}^* f(z)(y)$, by Theorem 4.6, we have $\|x\| \geq \|y\|$. One can show that, for any $z = (z_1, z_2) \in \mathbb{B}(\bar{z}, r)$, the following system of equations has solutions for $y_1, y_2$.

$$\begin{cases} y_1^2 + y_2^2 = 1, \\ y_1 z_1 z_2 = y_2 \frac{z_1^2 - z_2^2}{2}. \end{cases}$$

This implies that, for every $\eta > 0$, we have

$$\left\{\|x\|: x \in \widehat{D}^* f(z)(y), z \in \mathbb{B}(\bar{z}, \eta), f(z) \in \mathbb{B}(\bar{w}, \eta), \|y\| = 1, y_1 z_1 z_2 = y_2 \frac{z_1^2 - z_2^2}{2}\right\} \neq \emptyset.$$

Then, we calculate

$$\hat{\alpha}(f, \bar{z}, \bar{w}) = \sup_{\eta > 0} \inf \{\|x\|: x \in \widehat{D}^* f(z)(y), z \in \mathbb{B}(\bar{z}, \eta), f(z) \in \mathbb{B}(\bar{w}, \eta), \|y\| = 1\}$$

$$= \sup_{\eta > 0} \inf \left\{\|x\|: x \in \widehat{D}^* f(z)(y), z \in \mathbb{B}(\bar{z}, \eta), f(z) \in \mathbb{B}(\bar{w}, \eta), \|y\| = 1, y_1 z_1 z_2 = y_2 \frac{z_1^2 - z_2^2}{2}\right\}$$

$$= \sup_{\eta > 0} \inf \left\{\|y\|: x \in \widehat{D}^* f(z)(y), z \in \mathbb{B}(\bar{z}, \eta), f(z) \in \mathbb{B}(\bar{w}, \eta), \|y\| = 1, y_1 z_1 z_2 = y_2 \frac{z_1^2 - z_2^2}{2}\right\}$$

$$= \sup_{\eta > 0} \{1\} = 1. \qquad \square$$

Notice that in Theorem 3.3, for any $z = (z_1, z_2) \in \mathbb{R}^2 \setminus \{\theta\}$, the existence and solution of $\widehat{D}^* f(z)$ are proved by the Fréchet derivative $\nabla f(z)$ at $z$ and Theorem 1.38 in [16]. It is clearly to see that it is very difficult to solve $\widehat{D}^* f(z)$ without using the Fréchet derivative $\nabla f(z)$ at $z$. Next, we directly prove the results of $\widehat{D}^* f(z)$ in Part (b) of Theorem 3.3.

**Proposition 3.6.** *Let* $z = (z_1, z_2) \in \mathbb{R}^2 \setminus \{\theta\}$ *and* $x = (x_1, x_2), y = (y_1, y_2) \in \mathbb{R}^2$. *If* $x \in \widehat{D}^* f(z)(y)$, *then*

$$x_1 = y_1 \left(\frac{-(z_1^2 - z_2^2) z_1}{(z_1^2 + z_2^2)\sqrt{z_1^2 + z_2^2}} + \frac{2 z_1}{\sqrt{z_1^2 + z_2^2}}\right) + y_2 \left(\frac{-2 z_1 z_2 z_1}{(z_1^2 + z_2^2)\sqrt{z_1^2 + z_2^2}} + \frac{2 z_2}{\sqrt{z_1^2 + z_2^2}}\right)$$

$$= y_1 \frac{(z_1^2 + 3 z_2^2) z_1}{(z_1^2 + z_2^2)\sqrt{z_1^2 + z_2^2}} + y_2 \frac{2 z_2^2 z_2}{(z_1^2 + z_2^2)\sqrt{z_1^2 + z_2^2}},$$

and
$$x_2 = y_1 \left( \frac{-(z_1^2-z_2^2)z_2}{(z_1^2+z_2^2)\sqrt{z_1^2+z_2^2}} + \frac{-2z_2}{\sqrt{z_1^2+z_2^2}} \right) + y_2 \left( \frac{-2z_1 z_2 z_2}{(z_1^2+z_2^2)\sqrt{z_1^2+z_2^2}} + \frac{2z_1}{\sqrt{z_1^2+z_2^2}} \right)$$

$$= y_1 \frac{(-3z_1^2-z_2^2)z_2}{(z_1^2+z_2^2)\sqrt{z_1^2+z_2^2}} + y_2 \frac{2z_1^2 z_1}{(z_1^2+z_2^2)\sqrt{z_1^2+z_2^2}}.$$

*Proof.* Let $y = (y_1, y_2) \in \mathbb{R}^2$ and $x = (x_1, x_2) \in \mathbb{R}^2$, in order to check if $x \in \widehat{D}^* f(z)(y)$ or not, we calculate the following limits.

$$\limsup_{u \to z} \frac{\langle x,\ u-z \rangle - \langle y,\ f(u)-f(z) \rangle}{\|u-z\| + \|f(u)-f(z)\|}$$

$$= \limsup_{u \to z} \frac{\left\langle (x_1,x_2),\ (u_1,u_2)-(z_1,z_2) \right\rangle - \left\langle (y_1,y_2),\ \left( \frac{u_1^2-u_2^2}{\sqrt{u_1^2+u_2^2}}, \frac{2u_1 u_2}{\sqrt{u_1^2+u_2^2}} \right) - \left( \frac{z_1^2-z_2^2}{\sqrt{z_1^2+z_2^2}}, \frac{2z_1 z_2}{\sqrt{z_1^2+z_2^2}} \right) \right\rangle}{\|u-z\| + \|f(u)-f(z)\|}$$

$$= \limsup_{u \to z} \frac{x_1(u_1-z_1)+x_2(u_2-z_2) - y_1 \left( \frac{u_1^2-u_2^2}{\sqrt{u_1^2+u_2^2}} - \frac{z_1^2-z_2^2}{\sqrt{z_1^2+z_2^2}} \right) - y_2 \left( \frac{2u_1 u_2}{\sqrt{u_1^2+u_2^2}} - \frac{2z_1 z_2}{\sqrt{z_1^2+z_2^2}} \right)}{\|u-z\| + \|f(u)-f(z)\|}$$

$$= \limsup_{u \to z} \frac{x_1(u_1-z_1)+x_2(u_2-z_2) - y_1 \left( \frac{u_1^2-u_2^2}{\sqrt{u_1^2+u_2^2}} - \frac{z_1^2-z_2^2}{\sqrt{z_1^2+z_2^2}} \right) - y_2 \left( \frac{2u_1 u_2}{\sqrt{u_1^2+u_2^2}} - \frac{2z_1 z_2}{\sqrt{z_1^2+z_2^2}} \right)}{\sqrt{(u_1-z_1)^2+(u_2-z_2)^2} + \sqrt{\left( \frac{u_1^2-u_2^2}{\sqrt{u_1^2+u_2^2}} - \frac{z_1^2-z_2^2}{\sqrt{z_1^2+z_2^2}} \right)^2 + \left( \frac{2u_1 u_2}{\sqrt{u_1^2+u_2^2}} - \frac{2z_1 z_2}{\sqrt{z_1^2+z_2^2}} \right)^2}}. \quad (3.4)$$

In the limit (3.4), we consider some special directions for $u \to z$.

(D1) We take a special direction in the limit (3.4) as $u = (z_1, z_2 + s)$ with $s \downarrow 0$. We have

$$\limsup_{u \to z} \frac{\langle x,\ u-z \rangle - \langle y,\ f(u)-f(z) \rangle}{\|u-z\| + \|f(u)-f(z)\|}$$

$$\geq \lim_{u = (z_1, z_2+s)\ \text{with}\ s \downarrow 0} \frac{\langle x,\ u-z \rangle - \langle y,\ f(u)-f(z) \rangle}{\|u-z\| + \|f(u)-f(z)\|}$$

$$= \lim_{s \downarrow 0} \frac{s x_2 - y_1 \left( \frac{z_1^2-(z_2+s)^2}{\sqrt{z_1^2+(z_2+s)^2}} - \frac{z_1^2-z_2^2}{\sqrt{z_1^2+z_2^2}} \right) - y_2 \left( \frac{2z_1(z_2+s)}{\sqrt{z_1^2+(z_2+s)^2}} - \frac{2z_1 z_2}{\sqrt{z_1^2+z_2^2}} \right)}{\sqrt{(z_1-z_1)^2+(z_2+s-z_2)^2} + \sqrt{\left( \frac{z_1^2-(z_2+s)^2}{\sqrt{z_1^2+(z_2+s)^2}} - \frac{z_1^2-z_2^2}{\sqrt{z_1^2+z_2^2}} \right)^2 + \left( \frac{2z_1(z_2+s)}{\sqrt{z_1^2+(z_2+s)^2}} - \frac{2z_1 z_2}{\sqrt{z_1^2+z_2^2}} \right)^2}}$$

$$= \lim_{s \downarrow 0} \frac{s x_2 - y_1 \left( \frac{z_1^2-(z_2+s)^2}{\sqrt{z_1^2+(z_2+s)^2}} - \frac{z_1^2-z_2^2}{\sqrt{z_1^2+z_2^2}} \right) - y_2 \left( \frac{2z_1(z_2+s)}{\sqrt{z_1^2+(z_2+s)^2}} - \frac{2z_1 z_2}{\sqrt{z_1^2+z_2^2}} \right)}{s + \sqrt{\left( \frac{z_1^2-(z_2+s)^2}{\sqrt{z_1^2+(z_2+s)^2}} - \frac{z_1^2-z_2^2}{\sqrt{z_1^2+z_2^2}} \right)^2 + \left( \frac{2z_1(z_2+s)}{\sqrt{z_1^2+(z_2+s)^2}} - \frac{2z_1 z_2}{\sqrt{z_1^2+z_2^2}} \right)^2}}. \quad (3.5)$$

We calculate the two big terms in (3.5).

$$\frac{z_1^2-(z_2+s)^2}{\sqrt{z_1^2+(z_2+s)^2}}-\frac{z_1^2-z_2^2}{\sqrt{z_1^2+z_2^2}}$$

$$=\frac{(z_1^2-z_2^2)-2sz_2-s^2}{\sqrt{(z_1^2+z_2^2)+2sz_2+s^2}}-\frac{z_1^2-z_2^2}{\sqrt{z_1^2+z_2^2}}$$

$$=\frac{z_1^2-z_2^2}{\sqrt{(z_1^2+z_2^2)+2sz_2+s^2}}-\frac{z_1^2-z_2^2}{\sqrt{z_1^2+z_2^2}}+\frac{-2sz_2-s^2}{\sqrt{(z_1^2+z_2^2)+2sz_2+s^2}}$$

$$=\frac{(z_1^2-z_2^2)\left(\sqrt{z_1^2+z_2^2}-\sqrt{(z_1^2+z_2^2)+2sz_2+s^2}\right)}{\sqrt{(z_1^2+z_2^2)+2sz_2+s^2}\sqrt{z_1^2+z_2^2}}+\frac{-2sz_2-s^2}{\sqrt{(z_1^2+z_2^2)+2sz_2+s^2}}$$

$$=\frac{(z_1^2-z_2^2)\left(\left(\sqrt{z_1^2+z_2^2}\right)^2-\left(\sqrt{(z_1^2+z_2^2)+2sz_2+s^2}\right)^2\right)}{\sqrt{(z_1^2+z_2^2)+2sz_2+s^2}\sqrt{z_1^2+z_2^2}\left(\sqrt{z_1^2+z_2^2}+\sqrt{(z_1^2+z_2^2)+2sz_2+s^2}\right)}+\frac{-2sz_2-s^2}{\sqrt{(z_1^2+z_2^2)+2sz_2+s^2}}$$

$$=\frac{(z_1^2-z_2^2)\left(\left((z_1^2+z_2^2)\right)-\left((z_1^2+z_2^2)+2sz_2+s^2\right)\right)}{\sqrt{(z_1^2+z_2^2)+2sz_2+s^2}\sqrt{z_1^2+z_2^2}\left(\sqrt{z_1^2+z_2^2}+\sqrt{(z_1^2+z_2^2)+2sz_2+s^2}\right)}+\frac{-2sz_2-s^2}{\sqrt{(z_1^2+z_2^2)+2sz_2+s^2}}$$

$$=\frac{(z_1^2-z_2^2)(-2sz_2-s^2)}{\sqrt{(z_1^2+z_2^2)+2sz_2+s^2}\sqrt{z_1^2+z_2^2}\left(\sqrt{z_1^2+z_2^2}+\sqrt{(z_1^2+z_2^2)+2sz_2+s^2}\right)}+\frac{-2sz_2-s^2}{\sqrt{(z_1^2+z_2^2)+2sz_2+s^2}}$$

$$=s\left(\frac{-(z_1^2-z_2^2)(2z_2+s)}{\sqrt{(z_1^2+z_2^2)+2sz_2+s^2}\sqrt{z_1^2+z_2^2}\left(\sqrt{z_1^2+z_2^2}+\sqrt{(z_1^2+z_2^2)+2sz_2+s^2}\right)}+\frac{-2z_2-s}{\sqrt{(z_1^2+z_2^2)+2sz_2+s^2}}\right)$$

$$\equiv sA(s). \qquad (3.6)$$

Now we calculate other big term in (3.5).

$$\frac{2z_1(z_2+s)}{\sqrt{(z_1)^2+(z_2+s)^2}}-\frac{2z_1z_2}{\sqrt{z_1^2+z_2^2}}$$

$$=\frac{2z_1z_2+2sz_1}{\sqrt{(z_1^2+z_2^2)+2sz_2+s^2}}-\frac{2z_1z_2}{\sqrt{z_1^2+z_2^2}}$$

$$=\frac{2z_1z_2}{\sqrt{(z_1^2+z_2^2)+2sz_2+s^2}}-\frac{2z_1z_2}{\sqrt{z_1^2+z_2^2}}+\frac{2sz_1}{\sqrt{(z_1^2+z_2^2)+2sz_2+s^2}}$$

$$=\frac{2z_1z_2\left(\sqrt{z_1^2+z_2^2}-\sqrt{(z_1^2+z_2^2)+2s(z_1+z_2)+2s^2}\right)}{\sqrt{(z_1^2+z_2^2)+2sz_2+s^2}\sqrt{z_1^2+z_2^2}}+\frac{2sz_1}{\sqrt{(z_1^2+z_2^2)+2sz_2+s^2}}$$

$$= \frac{2z_1 z_2\Big((z_1^2+z_2^2)-\big((z_1^2+z_2^2)+2sz_2+s^2\big)\Big)}{\sqrt{(z_1^2+z_2^2)+2sz_2+s^2}\sqrt{z_1^2+z_2^2}\Big(\sqrt{z_1^2+z_2^2}+\sqrt{(z_1^2+z_2^2)+2sz_2+s^2}\Big)} + \frac{2sz_1}{\sqrt{(z_1^2+z_2^2)+2sz_2+s^2}}$$

$$= \frac{2z_1 z_2(-2sz_2-s^2))}{\sqrt{(z_1^2+z_2^2)+2sz_2+s^2}\sqrt{z_1^2+z_2^2}\Big(\sqrt{z_1^2+z_2^2}+\sqrt{(z_1^2+z_2^2)+2sz_2+s^2}\Big)} + \frac{2sz_1}{\sqrt{(z_1^2+z_2^2)+2sz_2+s^2}}$$

$$= s\left(\frac{-2z_1 z_2(2z_2+s)}{\sqrt{(z_1^2+z_2^2)+2sz_2+s^2}\sqrt{z_1^2+z_2^2}\Big(\sqrt{z_1^2+z_2^2}+\sqrt{(z_1^2+z_2^2)+2sz_2+s^2}\Big)} + \frac{2z_1}{\sqrt{(z_1^2+z_2^2)+2sz_2+s^2}}\right)$$

$$= sB(s). \tag{3.7}$$

Substituting (3.6) and (3.7) into the limit (3.5), we have

$$\limsup_{u\to z} \frac{\langle x,\ u-z\rangle - \langle y,\ f(u)-f(z)\rangle}{\|u-z\|+\|f(u)-f(z)\|}$$

$$\geq \lim_{s\downarrow 0} \frac{sx_2 - y_1 sA(s) - y_2 sB(s)}{s+\sqrt{(sA(s))^2+(sB(s))^2}}$$

$$= \lim_{s\downarrow 0} \frac{sx_2 - y_1 sA(s) - y_2 sB(s)}{s+ s\sqrt{(A(s))^2+(B(s))^2}}$$

$$= \lim_{s\downarrow 0} \frac{x_2 - y_1(A(s)) - y_2(B(s))}{1+\sqrt{(A(s))^2+(B(s))^2}}. \tag{3.8}$$

In the limit (3.8), by notations (3.6) and (3.7), we have

$$\lim_{s\downarrow 0} A(s) = \lim_{s\downarrow 0}\left(\frac{-(z_1^2-z_2^2)(2z_2+s)}{\sqrt{(z_1^2+z_2^2)+2sz_2+s^2}\sqrt{z_1^2+z_2^2}\Big(\sqrt{z_1^2+z_2^2}+\sqrt{(z_1^2+z_2^2)+2sz_2+s^2}\Big)} + \frac{-2z_2-s}{\sqrt{(z_1^2+z_2^2)+2sz_2+s^2}}\right)$$

$$= \frac{-(z_1^2-z_2^2)2z_2}{\sqrt{z_1^2+z_2^2}\sqrt{z_1^2+z_2^2}\Big(\sqrt{z_1^2+z_2^2}+\sqrt{z_1^2+z_2^2}\Big)} + \frac{-2z_2}{\sqrt{z_1^2+z_2^2}}$$

$$= \frac{-(z_1^2-z_2^2)z_2}{(z_1^2+z_2^2)\sqrt{z_1^2+z_2^2}} + \frac{-2z_2}{\sqrt{z_1^2+z_2^2}}. \tag{3.9}$$

We also have

$$\lim_{s\downarrow 0} B(s) = \lim_{s\downarrow 0}\left(\frac{-2z_1 z_2(2z_2+s)}{\sqrt{(z_1^2+z_2^2)+2sz_2+s^2}\sqrt{z_1^2+z_2^2}\Big(\sqrt{z_1^2+z_2^2}+\sqrt{(z_1^2+z_2^2)+2sz_2+s^2}\Big)} + \frac{2z_1}{\sqrt{(z_1^2+z_2^2)+2sz_2+s^2}}\right)$$

$$= \frac{-4z_1 z_2 z_2}{\sqrt{z_1^2+z_2^2}\sqrt{z_1^2+z_2^2}\Big(\sqrt{z_1^2+z_2^2}+\sqrt{z_1^2+z_2^2}\Big)} + \frac{2z_1}{\sqrt{z_1^2+z_2^2}}$$

$$= \frac{-2z_1 z_2 z_2}{(z_1^2+z_2^2)\sqrt{z_1^2+z_2^2}} + \frac{2z_1}{\sqrt{z_1^2+z_2^2}}. \qquad (3.10)$$

Substituting (3.9) and (3.10) into (3.8), we have

$$\limsup_{u \to z} \frac{\langle x,\ u-z \rangle - \langle y,\ f(u)-f(z) \rangle}{\|u-z\| + \|f(u)-f(z)\|}$$

$$\geq \lim_{s \downarrow 0} \frac{x_2 - y_1(E(s)) - y_2(F(s))}{1 + \sqrt{(E(s))^2 + (F(s))^2}}$$

$$= \frac{x_2 - y_1 \left( \frac{-(z_1^2 - z_2^2)z_2}{(z_1^2+z_2^2)\sqrt{z_1^2+z_2^2}} + \frac{-2z_2}{\sqrt{z_1^2+z_2^2}} \right) - y_2 \left( \frac{-2z_1 z_2 z_2}{(z_1^2+z_2^2)\sqrt{z_1^2+z_2^2}} + \frac{2z_1}{\sqrt{z_1^2+z_2^2}} \right)}{1 + \sqrt{\left( \frac{-(z_1^2 - z_2^2)z_2}{(z_1^2+z_2^2)\sqrt{z_1^2+z_2^2}} + \frac{-2z_2}{\sqrt{z_1^2+z_2^2}} \right)^2 + \left( \frac{-2z_1 z_2 z_2}{(z_1^2+z_2^2)\sqrt{z_1^2+z_2^2}} + \frac{2z_1}{\sqrt{z_1^2+z_2^2}} \right)^2}}.$$

This implies that

$$x_2 - y_1 \left( \frac{-(z_1^2 - z_2^2)z_2}{(z_1^2+z_2^2)\sqrt{z_1^2+z_2^2}} + \frac{-2z_2}{\sqrt{z_1^2+z_2^2}} \right) - y_2 \left( \frac{-2z_1 z_2 z_2}{(z_1^2+z_2^2)\sqrt{z_1^2+z_2^2}} + \frac{2z_1}{\sqrt{z_1^2+z_2^2}} \right) > 0 \Rightarrow x \notin \widehat{D}^* f(z)(y).$$

This reduces to that

(i) $x \in \widehat{D}^* f(z)(y) \Rightarrow x_2 - y_1 \left( \frac{-(z_1^2 - z_2^2)z_2}{(z_1^2+z_2^2)\sqrt{z_1^2+z_2^2}} + \frac{-2z_2}{\sqrt{z_1^2+z_2^2}} \right) - y_2 \left( \frac{-2z_1 z_2 z_2}{(z_1^2+z_2^2)\sqrt{z_1^2+z_2^2}} + \frac{2z_1}{\sqrt{z_1^2+z_2^2}} \right) \leq 0.$

(D2) We take a new special direction in the limit (3.4) $u = (z_1, z_2 - s)$ with $s \downarrow 0$ That is $s > 0$. With respect to this direction of the limit in (3.4), similarly, to the proof of (ii), we can show that

(j) $x \in \widehat{D}^* f(z)(y) \Rightarrow -x_2 - y_1 \left( \frac{(z_1^2 - z_2^2)z_2}{(z_1^2+z_2^2)\sqrt{z_1^2+z_2^2}} + \frac{2z_2}{\sqrt{z_1^2+z_2^2}} \right) - y_2 \left( \frac{2z_1 z_2 z_2}{(z_1^2+z_2^2)\sqrt{z_1^2+z_2^2}} + \frac{-2z_1}{\sqrt{z_1^2+z_2^2}} \right) \leq 0.$

Combining (i) and (j), we obtain that

(I) $\quad x \in \widehat{D}^* f(z)(y) \Rightarrow x_2 = y_1 \left( \frac{-(z_1^2 - z_2^2)z_2}{(z_1^2+z_2^2)\sqrt{z_1^2+z_2^2}} + \frac{-2z_2}{\sqrt{z_1^2+z_2^2}} \right) + y_2 \left( \frac{-2z_1 z_2 z_2}{(z_1^2+z_2^2)\sqrt{z_1^2+z_2^2}} + \frac{2z_1}{\sqrt{z_1^2+z_2^2}} \right).$

We take new special directions in the limit (3.4):

**(D3)** $u = (z_1 + s, z_2)$ with $s \downarrow 0$ That is $s > 0$ and **(D4)** $u = (z_1 - s, z_2)$ with $s \downarrow 0$ That is $s > 0$. Similarly, to the proof of (I), we can show that

(II) $\quad x \in \widehat{D}^* f(z)(y) \Rightarrow x_1 = y_1 \left( \frac{-(z_1^2 - z_2^2)z_1}{(z_1^2+z_2^2)\sqrt{z_1^2+z_2^2}} + \frac{2z_1}{\sqrt{z_1^2+z_2^2}} \right) + y_2 \left( \frac{-2z_1 z_2 z_1}{(z_1^2+z_2^2)\sqrt{z_1^2+z_2^2}} + \frac{2z_2}{\sqrt{z_1^2+z_2^2}} \right).$

By (I) and (II), we have that $x \in \widehat{D}^*f(z)(y)$ implies that

$$x_1 = y_1\left(\frac{-(z_1^2-z_2^2)z_1}{(z_1^2+z_2^2)\sqrt{z_1^2+z_2^2}} + \frac{2z_1}{\sqrt{z_1^2+z_2^2}}\right) + y_2\left(\frac{-2z_1z_2z_1}{(z_1^2+z_2^2)\sqrt{z_1^2+z_2^2}} + \frac{2z_2}{\sqrt{z_1^2+z_2^2}}\right),$$

$$x_2 = y_1\left(\frac{-(z_1^2-z_2^2)z_2}{(z_1^2+z_2^2)\sqrt{z_1^2+z_2^2}} + \frac{-2z_2}{\sqrt{z_1^2+z_2^2}}\right) + y_2\left(\frac{-2z_1z_2z_2}{(z_1^2+z_2^2)\sqrt{z_1^2+z_2^2}} + \frac{2z_1}{\sqrt{z_1^2+z_2^2}}\right). \qquad \square$$

## 4. Fréchet and Mordukhovich Differentiability of Mappings in $\mathbb{R}^4$

In this section, we extend the mapping studied in the previous section from $\mathbb{R}^2$ to $\mathbb{R}^4$. We first consider a mapping $g: \mathbb{R}^4 \to \mathbb{R}^4$ by

$$g((x_1, x_2, x_3, x_4)) = \left(\frac{x_1^2 - x_2^2}{\sqrt{x_1^2+x_2^2}}, \frac{2x_1x_2}{\sqrt{x_1^2+x_2^2}}, \frac{x_3^2 - x_4^2}{\sqrt{x_3^2+x_4^2}}, \frac{2x_3x_4}{\sqrt{x_3^2+x_4^2}}\right), \text{ for } (x_1, x_2, x_3, x_4) \in \mathbb{R}^4. \qquad (4.1)$$

Here, for $i = 1, 3$, if $x_i = x_{i+1} = 0$, then $\left(\frac{x_i^2 - x_{i+1}^2}{\sqrt{x_i^2+x_{i+1}^2}}, \frac{2x_ix_{i+1}}{\sqrt{x_i^2+x_{i+1}^2}}\right)$ is defined to be $(0, 0)$.

Similar, to the norm preserving property of $f$ studied in the previous section, the mapping $g: \mathbb{R}^4 \to \mathbb{R}^4$ is also norm preserving satisfying,

$$\|g(x)\| = \|x\|, \text{ for any } x \in \mathbb{R}^4.$$

**Lemma 4.1.** *Let $z = (z_1, z_2, z_3, z_4) \in \mathbb{R}^4$. If $(z_1^2 + z_2^2)(z_3^2 + z_4^2) = 0$, then $g$ is not Fréchet differentiable at $z$.*

*Proof.* The proof of this lemma is similar to the proof of Lemma 3.1. But it is more complicated. So we prove it here. Let $z = (z_1, z_2, z_3, z_4) \in \mathbb{R}^4$ with $(z_1^2 + z_2^2)(z_3^2 + z_4^2) = 0$.

Suppose that $z_1^2 + z_2^2 = 0$. Assume, by the way of contradiction, that $g$ is Fréchet differentiable at $z$. Then $\nabla g(z): \mathbb{R}^4 \to \mathbb{R}^4$ must be precisely represented by a real $4 \times 4$ matrix as below.

$$\nabla g(z) = \begin{pmatrix} a_{11} & a_{12} & a_{13} & a_{14} \\ a_{21} & a_{22} & a_{23} & a_{24} \\ a_{31} & a_{32} & a_{33} & a_{34} \\ a_{41} & a_{42} & a_{43} & a_{44} \end{pmatrix} \equiv A.$$

By the assumption that $z_1^2 + z_2^2 = 0$, that is, $z = (0, 0, z_3, z_4)$, we have

$$\lim_{u=(u_1,u_2,u_3,u_4) \to \theta} \frac{g(u) - g(z) - \nabla g(z)(u-z)}{\|u-z\|}$$

$$= \lim_{u \to \theta} \frac{\left(\frac{u_1^2-u_2^2}{\sqrt{u_1^2+u_2^2}} - 0, \frac{2u_1u_2}{\sqrt{u_1^2+u_2^2}} - 0, \frac{u_3^2-u_4^2}{\sqrt{u_3^2+u_4^2}} - \frac{z_3^2-z_4^2}{\sqrt{z_3^2+z_4^2}}, \frac{2u_3u_4}{\sqrt{u_3^2+u_4^2}} - \frac{2z_3z_4}{\sqrt{z_3^2+z_4^2}}\right) - (u-z)A}{\|u-z\|}. \qquad (4.2)$$

Assume $a_{11} \neq 0$. In this case, we take a special direction in the limit (4.2) for $u \to \theta$ by $u_1 = -a_{11}t$ with $t \downarrow 0$, $u_2 = 0$, $u_3 = z_3$ and $u_4 = z_4$. By the assumption that $\nabla g(z) = A$, it follows that

$$\lim_{u=(u_1,u_2,u_3,u_4) \to \theta} \frac{g(u)-g(z)-\nabla g(z)(u-z)}{\|u-z\|}$$

$$= \lim_{t \downarrow 0} \frac{\left(\frac{a_{11}^2 t^2}{\sqrt{a_{11}^2 t^2}},\ 0,\ 0,\ 0\right) - (-a_{11}t,\ 0,\ 0,\ 0)A}{\|(-a_{11}t,0,0,0)\|}$$

$$= \lim_{t \downarrow 0} \frac{(|a_{11}|t,\ 0,\ 0,\ 0)+(a_{11}^2 t,\ a_{11}a_{12}t,\ a_{11}a_{13}t,\ a_{11}a_{14}t)}{|a_{11}|t}$$

$$= (1+|a_{11}|,\ (\operatorname{sign} a_{11})a_{12},\ (\operatorname{sign} a_{11})a_{13},\ (\operatorname{sign} a_{11})a_{14})$$

$$\neq \theta.$$

This contradicts the assumption that $\nabla g(z) = A$ (which implies that $a_{11} = 0$).

Next, suppose $a_{11} = 0$. In this case, we take a special direction in the limit (4.2) for $u \to \theta$ by $u_1 \downarrow 0$, $u_2 = 0$, $u_3 = z_3$ and $u_4 = z_4$. It follows that

$$\lim_{u=(u_1,u_2,u_3,u_4) \to \theta} \frac{g(u)-g(z)-\nabla g(z)(u-z)}{\|u-z\|}$$

$$= \lim_{u_1 \downarrow 0} \frac{\left(\frac{u_1^2}{\sqrt{u_1^2}},\ 0,\ 0,\ 0\right) - (u_1,\ 0,\ 0,\ 0)A}{\|(u_1,0,0,0)\|}$$

$$= \lim_{u_1 \downarrow 0} \frac{(u_1,\ 0,\ 0,\ 0)+(0,\ u_1 a_{12},\ u_1 a_{13},\ u_1 a_{14})}{u_1}$$

$$= (1,\ a_{12},\ u_1 a_{13},\ u_1 a_{14})$$

$$\neq \theta.$$

This contradicts the assumption that $\nabla g(z) = A$.

Next, we suppose that $z_3^2 + z_4^2 = 0$. In this case, $z = (z_1, z_2, 0, 0)$, we have

$$\lim_{u=(u_1,u_2,u_3,u_4) \to \theta} \frac{g(u)-g(z)-\nabla g(z)(u-z)}{\|u-z\|}$$

$$= \lim_{u \to \theta} \frac{\left(\frac{u_1^2-u_2^2}{\sqrt{u_1^2+u_2^2}} - \frac{z_1^2-z_2^2}{\sqrt{z_1^2+z_2^2}},\ \frac{2u_1 u_2}{\sqrt{u_1^2+u_2^2}} - \frac{2z_1 z_2}{\sqrt{z_1^2+z_2^2}},\ \frac{u_3^2-u_4^2}{\sqrt{u_3^2+u_4^2}} - 0,\ \frac{2u_3 u_4}{\sqrt{u_3^2+u_4^2}} - 0\right) - (u-z)A}{\|u-z\|}. \qquad (4.3)$$

Assume $a_{33} \neq 0$. In this case, we take a special direction in the limit (4.2) for $u \to \theta$ by $u_3 = -a_{33}t$ with $t \downarrow 0$, $u_4 = 0$, $u_2 = z_2$ and $u_1 = z_1$. By the assumption that $\nabla g(z) = A$, it follows that

$$\lim_{u=(u_1,u_2,u_3,u_4) \to \theta} \frac{g(u)-g(z)-\nabla g(z)(u-z)}{\|u-z\|}$$

$$= \lim_{t \downarrow 0} \frac{\left(0,\ 0,\ \frac{a_{33}^2 t^2}{\sqrt{a_{33}^2 t^2}},\ 0\right) - (0,\ 0,\ -a_{33}t,\ 0)A}{\|(0, 0, -a_{33}t, 0)\|}$$

$$= \lim_{t \downarrow 0} \frac{(0,\ 0,\ |a_{33}|t,\ 0) + (a_{31}a_{33}t,\ a_{32}a_{33}t,\ a_{33}^2 t,\ a_{34}a_{33}t)}{|a_{33}|t}$$

$$= ((\operatorname{sign} a_{33})a_{31}, (\operatorname{sign} a_{32})a_{32}, 1 + |a_{33}|, (\operatorname{sign} a_{11})a_{34})$$

$$\neq \theta.$$

This contradicts the assumption that $\nabla g(z) = A$.

Next, suppose $a_{33} = 0$. In this case, we take a special direction in the limit (4.3) for $u \to \theta$ by $u_3 \downarrow 0$, $u_4 = 0$, $u_2 = z_2$ and $u_1 = z_1$. It follows that

$$\lim_{u = (u_1, u_2, u_3, u_4) \to \theta} \frac{g(u) - g(z) - \nabla g(z)(u - z)}{\|u - z\|}$$

$$= \lim_{u_3 \downarrow 0} \frac{\left(0,\ 0,\ \frac{u_3^2}{\sqrt{u_3^2}},\ 0\right) - (0,\ 0,\ u_3,\ 0)A}{\|(0, 0, u_3, 0)\|}$$

$$= \lim_{u_3 \downarrow 0} \frac{(0,\ 0,\ u_3,\ 0) + (u_3 a_{31},\ u_3 a_{32},\ u_3 0,\ u_3 a_{34})}{u_3}$$

$$= (a_{31}, a_{32}, 1, a_{34})$$

$$\neq \theta.$$

This contradicts the assumption that $\nabla g(z) = A$. It completes the proof. □

**Corollary 4.2.** $g$ is not Fréchet differentiable at $\theta$.

*Proof.* In Lemma 4.1, let $z = (z_1, z_2, z_3, z_4) \in \mathbb{R}^4$ with $(z_1^2 + z_2^2) = 0$ and $(z_3^2 + z_4^2) = 0$. □

**Proposition 4.3.** *The mapping $g \colon \mathbb{R}^4 \to \mathbb{R}^4$ defined by (4.1) satisfies that*

(a) $\widehat{D}^* g(\theta)(\theta) = \{\theta\}$;
(b) $\widehat{D}^* g(\theta)(y) = \emptyset$, *for any $y \in \mathbb{R}^4 \setminus \{\theta\}$.*

*Proof.* The proof of this proposition is similar to the proof of Propositions 3.2. It is omitted here. □

**Theorem 4.4.** *Let $z = (z_1, z_2, z_3, z_4) \in \mathbb{R}^4 \setminus \{\theta\}$ with $(z_1^2 + z_2^2)(z_3^2 + z_4^2) \neq 0$. Then we have that*

(a) $g$ *is Fréchet differentiable at $z$ and*

$$\nabla g(z) = \begin{pmatrix} \frac{(z_1^2+3z_2^2)z_1}{(z_1^2+z_2^2)\sqrt{z_1^2+z_2^2}} & \frac{2z_2^2 z_2}{(z_1^2+z_2^2)\sqrt{z_1^2+z_2^2}} & 0 & 0 \\ \frac{-(3z_1^2+z_2^2)z_2}{(z_1^2+z_2^2)\sqrt{z_1^2+z_2^2}} & \frac{2z_1^2 z_1}{(z_1^2+z_2^2)\sqrt{z_1^2+z_2^2}} & 0 & 0 \\ 0 & 0 & \frac{(z_3^2+3z_4^2)z_3}{(z_3^2+z_4^2)\sqrt{z_3^2+z_4^2}} & \frac{2z_4^2 z_4}{(z_3^2+z_4^2)\sqrt{z_3^2+z_4^2}} \\ 0 & 0 & \frac{-(3z_3^2+z_4^2)z_4}{(z_3^2+z_4^2)\sqrt{z_3^2+z_4^2}} & \frac{2z_3^2 z_3}{(z_3^2+z_4^2)\sqrt{z_3^2+z_4^2}} \end{pmatrix}.$$

(b) $g$ is Mordukhovich differentiable at $z$ and

$$\widehat{D}^*g(z) = \begin{pmatrix} \frac{(z_1^2+3z_2^2)z_1}{(z_1^2+z_2^2)\sqrt{z_1^2+z_2^2}} & \frac{-(3z_1^2+z_2^2)z_2}{(z_1^2+z_2^2)\sqrt{z_1^2+z_2^2}} & 0 & 0 \\ \frac{2z_2^2 z_2}{(z_1^2+z_2^2)\sqrt{z_1^2+z_2^2}} & \frac{2z_1^2 z_1}{(z_1^2+z_2^2)\sqrt{z_1^2+z_2^2}} & 0 & 0 \\ 0 & 0 & \frac{(z_3^2+3z_4^2)z_3}{(z_3^2+z_4^2)\sqrt{z_3^2+z_4^2}} & \frac{-(3z_3^2+z_4^2)z_4}{(z_3^2+z_4^2)\sqrt{z_3^2+z_4^2}} \\ 0 & 0 & \frac{2z_4^2 z_4}{(z_3^2+z_4^2)\sqrt{z_3^2+z_4^2}} & \frac{2z_3^2 z_3}{(z_3^2+z_4^2)\sqrt{z_3^2+z_4^2}} \end{pmatrix}.$$

(More precisely speaking) Let $x = (x_1, x_2, x_3, x_4), y = (y_1, y_2, y_3, y_4) \in \mathbb{R}^4$. If $x = \widehat{D}^*g(z)(y)$, then

$$x_1 = y_1 \frac{(z_1^2+3z_2^2)z_1}{(z_1^2+z_2^2)\sqrt{z_1^2+z_2^2}} + y_2 \frac{2z_2^2 z_2}{(z_1^2+z_2^2)\sqrt{z_1^2+z_2^2}},$$

$$x_2 = y_1 \frac{(-3z_1^2-z_2^2)z_2}{(z_1^2+z_2^2)\sqrt{z_1^2+z_2^2}} + y_2 \frac{2z_1^2 z_1}{(z_1^2+z_2^2)\sqrt{z_1^2+z_2^2}},$$

$$x_3 = y_3 \frac{(z_3^2+3z_4^2)z_3}{(z_3^2+z_4^2)\sqrt{z_3^2+z_4^2}} + y_4 \frac{2z_4^2 z_4}{(z_3^2+z_4^2)\sqrt{z_3^2+z_4^2}},$$

$$x_4 = y_3 \frac{-(3z_3^2+z_4^2)z_4}{(z_3^2+z_4^2)\sqrt{z_3^2+z_4^2}} + y_4 \frac{2z_3^2 z_3}{(z_3^2+z_4^2)\sqrt{z_3^2+z_4^2}}.$$

*Proof.* The proof of this theorem is similar to the proof of Theorem 3.3. It is omitted here. □

**Theorem 4.5.** *Let $z = (z_1, z_2, z_3, z_4) \in \mathbb{R}^4 \setminus \{0\}$ with $(z_1^2 + z_2^2)(z_3^2 + z_4^2) = 0$. Let $x = (x_1, x_2, x_3, x_4)$ and $y = (y_1, y_2, y_3, y_4) \in \mathbb{R}^4$.*

(a) *Suppose that $z_1^2 + z_2^2 = 0$ ($z_3^2 + z_4^2 \neq 0$).*

(i) *If $y_1^2 + y_2^2 = 0$ and $x \in \widehat{D}^*g(z)(y)$, then $x$ and $y$ satisfy that*

$$x_1 = y_1 = 0,$$
$$x_2 = y_2 = 0,$$

$$x_3 = y_3 \frac{(z_3^2+3z_4^2)z_3}{(z_3^2+z_4^2)\sqrt{z_3^2+z_4^2}} + y_4 \frac{2z_4^2 z_4}{(z_3^2+z_4^2)\sqrt{z_3^2+z_4^2}},$$

$$x_4 = y_3 \frac{-(3z_3^2+z_4^2)z_4}{(z_3^2+z_4^2)\sqrt{z_3^2+z_4^2}} + y_4 \frac{2z_3^2 z_3}{(z_3^2+z_4^2)\sqrt{z_3^2+z_4^2}}.$$

(ii) If $y_1^2 + y_2^2 \neq 0$, then $\widehat{D}^* g(z)(y) = \emptyset$.

(b) *Suppose that* $z_3^2 + z_4^2 = 0$ ($z_1^2 + z_2^2 \neq 0$).

(i) If $y_3^2 + y_4^2 = 0$ and $x \in \widehat{D}^* g(z)(y)$, then $x$ and $y$ satisfy that

$$x_1 = y_1 \frac{(z_1^2+3z_2^2)z_1}{(z_1^2+z_2^2)\sqrt{z_1^2+z_2^2}} + y_2 \frac{2z_2^2 z_2}{(z_1^2+z_2^2)\sqrt{z_1^2+z_2^2}},$$

$$x_2 = y_1 \frac{(-3z_1^2-z_2^2)z_2}{(z_1^2+z_2^2)\sqrt{z_1^2+z_2^2}} + y_2 \frac{2z_1^2 z_1}{(z_1^2+z_2^2)\sqrt{z_1^2+z_2^2}}$$

$$x_3 = y_3 = 0,$$
$$x_4 = y_4 = 0,$$

(ii) If $y_3^2 + y_4^2 \neq 0$, then $\widehat{D}^* g(z)(y) = \emptyset$.

*Proof.* Proof of (a). Let $y = (y_1, y_2, y_3, y_4) \in \mathbb{R}^4$ and $x = (x_1, x_2, x_3, x_4) \in \mathbb{R}^4$, in order to check if $x \in \widehat{D}^* f(z)(y)$ or not, we take a special direction in the limit (4.2) as $u = (u_1, u_2, z_3, z_4)$ with $u \to z$, which is equivalent to $(u_1, u_2) \to (z_1, z_2)$. Let $f: \mathbb{R}^2 \to \mathbb{R}^2$ be defined by (1.1) in [15]. By definition (5.1), we have

$$g(u) - g(z) = \left(\frac{u_1^2-u_2^2}{\sqrt{u_1^2+u_2^2}}, \frac{2u_1 u_2}{\sqrt{u_1^2+u_2^2}}, \frac{z_3^2-z_4^2}{\sqrt{z_3^2+z_4^2}}, \frac{2z_3 z_4}{\sqrt{z_3^2+z_4^2}}\right) - \left(\frac{z_1^2-z_2^2}{\sqrt{z_1^2+z_2^2}}, \frac{2z_1 z_2}{\sqrt{z_1^2+z_2^2}}, \frac{z_3^2-z_4^2}{\sqrt{z_3^2+z_4^2}}, \frac{2z_3 z_4}{\sqrt{z_3^2+z_4^2}}\right).$$

This implies that

$$\langle y, g(u) - g(z) \rangle = \langle (y_1, y_2), \left(\frac{u_1^2-u_2^2}{\sqrt{u_1^2+u_2^2}}, \frac{2u_1 u_2}{\sqrt{u_1^2+u_2^2}}\right) - \left(\frac{z_1^2-z_2^2}{\sqrt{z_1^2+z_2^2}}, \frac{2z_1 z_2}{\sqrt{z_1^2+z_2^2}}\right) \rangle$$

$$= \langle (y_1, y_2), f(u_1, u_2) - f(z_1, z_2) \rangle.$$

With respect to this special limit, the limit (4.2) for $g$ in $\mathbb{R}^4$ can be converted to a limit in $\mathbb{R}^2$ with respect to the mapping $f: \mathbb{R}^2 \to \mathbb{R}^2$ be defined by (1.1).

$$\limsup_{u \to z} \frac{\langle x, u-z \rangle - \langle y, g(u)-g(z) \rangle}{\|u-z\| + \|g(u)-g(z)\|}$$

$$\geq \lim_{u=(u_1,u_2,z_3,z_4) \to z} \frac{\langle x, u-z \rangle - \langle y, g(u)-g(z) \rangle}{\|u-z\|+\|g(u)-g(z)\|}$$

$$= \lim_{(u_1,u_2)\to(z_1,z_2)} \frac{\langle(x_1,x_2), (u_1,u_2)-(z_1,z_2)\rangle - \langle(y_1,y_2), f(u_1,u_2)-f(z_1,z_2)\rangle}{\|(u_1,u_2)-(z_1,z_2)\|+\|f(u_1,u_2)-f(z_1,z_2)\|}.$$

Then, with respect to the case $z_1^2 + z_2^2 = 0$ (This implies that $(z_1, z_2)$ is the origin in $\mathbb{R}^2$), Part (i) can be proved by Propositions 3.2. Part (ii) can be similarly proved. □

**Theorem 4.6.** *Let* $z = (z_1, z_2, z_3, z_4) \in \mathbb{R}^4 \setminus \{\theta\}$ *with* $(z_1^2 + z_2^2)(z_3^2 + z_4^2) \neq 0$. *Let* $x = (x_1, x_2, x_3, x_4)$ *and* $y = (y_1, y_2, y_3, y_4) \in \mathbb{R}^4$ *If* $x \in \widehat{D}^*g(z)(y)$, *then $x$ and $y$ satisfy the following conditions.*

(i) $\|x\|^2 = \|y\|^2 + 12\left(y_1 \frac{z_1 z_2}{z_1^2+z_2^2} - y_2 \frac{z_1^2-z_2^2}{2(z_1^2+z_2^2)}\right)^2 + 12\left(y_3 \frac{z_3 z_4}{z_3^2+z_4^2} - y_4 \frac{z_3^2-z_4^2}{2(z_3^2+z_4^2)}\right)^2.$

(ii) $\|x\| \geq \|y\|.$

(iii) $y_1 z_1 z_2 = y_2 \frac{z_1^2-z_2^2}{2}$ *and* $y_3 z_3 z_4 = y_4 \frac{z_3^2-z_4^2}{2}$ $\Longrightarrow$ $\|x\| = \|y\|.$

*Proof.* Suppose that $x \in \widehat{D}^*g(z)(y)$. By Theorem 3.3 and by the proof of Proposition 3.4, we have

$$x_1^2 + x_2^2 + x_3^2 + x_4^2$$

$$= \left(y_1 \frac{(z_1^2+3z_2^2)z_1}{(z_1^2+z_2^2)\sqrt{z_1^2+z_2^2}} + y_2 \frac{2z_2^2 z_2}{(z_1^2+z_2^2)\sqrt{z_1^2+z_2^2}}\right)^2 + \left(y_1 \frac{(-3z_1^2-z_2^2)z_2}{(z_1^2+z_2^2)\sqrt{z_1^2+z_2^2}} + y_2 \frac{2z_1^2 z_1}{(z_1^2+z_2^2)\sqrt{z_1^2+z_2^2}}\right)^2$$

$$+ \left(y_3 \frac{(z_3^2+3z_4^2)z_3}{(z_3^2+z_4^2)\sqrt{z_3^2+z_4^2}} + y_4 \frac{2z_4^2 z_4}{(z_3^2+z_4^2)\sqrt{z_3^2+z_4^2}}\right)^2 + \left(y_3 \frac{-(3z_3^2+z_4^2)z_4}{(z_3^2+z_4^2)\sqrt{z_3^2+z_4^2}} + y_4 \frac{2z_3^2 z_3}{(z_3^2+z_4^2)\sqrt{z_3^2+z_4^2}}\right)^2$$

$$= y_1^2 + y_2^2 + 12\left(y_1 \frac{z_1 z_2}{z_1^2+z_2^2} - y_2 \frac{z_1^2-z_2^2}{2(z_1^2+z_2^2)}\right)^2 + y_3^2 + y_4^2 + 12\left(y_3 \frac{z_3 z_4}{z_3^2+z_4^2} - y_4 \frac{z_3^2-z_4^2}{2(z_3^2+z_4^2)}\right)^2. \quad \square$$

**Corollary 4.7.** *Let* $z = (z_1, z_2, z_3, z_4) \in \mathbb{R}^4 \setminus \{\theta\}$ *with* $(z_1^2 + z_2^2)(z_3^2 + z_4^2) = 0$. *Let* $x = (x_1, x_2, x_3, x_4)$ *and* $y = (y_1, y_2, y_3, y_4) \in \mathbb{R}^4$.

(a) *Suppose that* $z_1^2 + z_2^2 = 0$ $(z_3^2 + z_4^2 \neq 0)$. *If* $y_1^2 + y_2^2 = 0$ *and* $x \in \widehat{D}^*g(z)(y)$, *then*

(i) $\|x\|^2 = \|y\|^2 + 12\left(y_3 \frac{z_3 z_4}{z_3^2+z_4^2} - y_4 \frac{z_3^2-z_4^2}{2(z_3^2+z_4^2)}\right)^2.$

(ii) $\|x\| \geq \|y\|.$

(iii) $y_3 z_3 z_4 = y_4 \frac{z_3^2-z_4^2}{2}$ $\Longrightarrow$ $\|x\| = \|y\|.$

(b) *Suppose that* $z_3^2 + z_4^2 = 0$ $(z_1^2 + z_2^2 \neq 0)$. *If* $y_3^2 + y_4^2 = 0$ *and* $x \in \widehat{D}^*g(z)(y)$, *then*

(i) $\|x\|^2 = \|y\|^2 + 12\left(y_1 \frac{z_1 z_2}{z_1^2+z_2^2} - y_2 \frac{z_1^2-z_2^2}{2(z_1^2+z_2^2)}\right)^2.$

(ii) $\|x\| \geq \|y\|$.

(iii) $y_1 z_1 z_2 = y_2 \frac{z_1^2 - z_2^2}{2} \implies \|x\| = \|y\|$.

*Proof.* Notice that *if* $y_1^2 + y_2^2 = 0$, then $\|y\|^2 = y_3^2 + y_4^2$. If $y_3^2 + y_4^2 = 0$, then $\|y\|^2 = y_1^2 + y_2^2$. Then, this corollary follows from Theorems 4.5 and 4.6. □

**Theorem 4.8.** *Let* $g \colon \mathbb{R}^4 \to \mathbb{R}^4$ *be defined by* (5.1). *Then*
$$\hat{\alpha}(g(\bar{z})) = 1, \text{for } \bar{z} = (\bar{z}_1, \bar{z}_2, \bar{z}_3, \bar{z}_4) \in \mathbb{R}^4 \setminus \{\theta\}.$$

*Proof.* For any $\bar{z}, \bar{w} \in \mathbb{R}^4$ with $g(\bar{z}) = \bar{w}$ and $\eta > 0$, let $\mathbb{B}(\bar{z}, \eta)$ and $\mathbb{B}(\bar{w}, \eta)$ denote the closed balls in $\mathbb{R}^4$ with radius $\eta$ and centered at point $\bar{z}$ and $\bar{w}$, respectively. By definition, we have

$$\hat{\alpha}(g(\bar{z})) = \sup_{\eta > 0} \inf\{\|x\| : x \in \widehat{D}^* g(z)(y), z \in \mathbb{B}(\bar{z}, \eta), g(z) \in \mathbb{B}(\bar{w}, \eta), \|y\| = 1\}$$

$$= \sup_{\eta > 0} \inf\{\|x\| : x \in \widehat{D}^* g(z)(y), z \in \mathbb{B}(\bar{z}, \eta), g(z) \in \mathbb{B}(\bar{w}, \eta), \|y\| = 1\}. \tag{4.4}$$

In case $\theta \in \mathbb{B}(\bar{z}, \eta)$ and $g(\theta) = \theta \in \mathbb{B}(\bar{w}, \eta)$, then by Proposition 4.4, we have that

$$\widehat{D}^* g(\theta)(y) = \emptyset, \text{ for any } y \in \mathbb{R}^4 \text{ with } \|y\| = 1. \tag{4.5}$$

Let $z = (z_1, z_2, z_3, z_4) \in \mathbb{R}^4 \setminus \{\theta\}$ with $(z_1^2 + z_2^2)(z_3^2 + z_4^2) = 0$. Suppose that $z \in \mathbb{B}(\bar{z}, \eta)$, for some $\eta > 0$. Let $y = (y_1, y_2, y_3, y_4) \in \mathbb{R}^4$ with $\|y\| = 1$. By Theorem 4.5, we have that

$$z_1^2 + z_2^2 = 0 \text{ and } y_1^2 + y_2^2 \neq 0 \implies \widehat{D}^* g(z)(y) = \emptyset. \tag{4.6}$$

And $\qquad\qquad z_3^2 + z_4^2 = 0 \text{ and } y_3^2 + y_4^2 \neq 0 \implies \widehat{D}^* g(z)(y) = \emptyset. \tag{4.7}$

From (4.5), (4.6) and (4.7), we see that in the limit set in (4.4), there are some point $z \in \mathbb{B}(\bar{z}, \eta)$ and $y \in \mathbb{R}^4$ with $\|y\| = 1$ such that $\widehat{D}^* g(z)(y) = \emptyset$. Since $g \colon \mathbb{R}^4 \to \mathbb{R}^4$ defined by (4.1) is continuous on $\mathbb{R}^4$, it yields that, for every $\eta > 0$,

$$\{z \in \mathbb{B}(\bar{z}, \eta) : g(z) \in \mathbb{B}(\bar{w}, \eta), \|y\| = 1\} \neq \emptyset.$$

By Theorems 4.4 and 4.5, this implies that

$$\{\|x\| : x \in \widehat{D}^* g(z)(y), z \in \mathbb{B}(\bar{z}, \eta), g(z) \in \mathbb{B}(\bar{w}, \eta), \|y\| = 1\} \neq \emptyset, \text{ for every } \eta > 0.$$

However, if $\widehat{D}^* g(z)(y) \neq \emptyset$, for $z \in \mathbb{B}(\bar{z}, \eta)$ and $y \in \mathbb{R}^4$ with $\|y\| = 1$, then, for any $x \in \widehat{D}^* g(z)(y)$, by Theorem 4.6 and Corollary 4.7, we have $\|x\| \geq \|y\|$. One can show that, for any $z = (z_1, z_2, z_3, z_4) \in \mathbb{B}(\bar{z}, r)$, the following system of equations has solutions for $y_1, y_2, y_3, y_4$.

$$\begin{cases} y_1^2 + y_2^2 + y_3^2 + y_4^2 = 1, \\ y_1 z_1 z_2 = y_2 \frac{z_1^2 - z_2^2}{2}, \\ y_3 z_3 z_4 = y_4 \frac{z_3^2 - z_4^2}{2}. \end{cases}$$

Meanwhile, following systems of equations have solutions for $y_1, y_2$ and for $y_3, y_4$, respectively:

$$\begin{cases} y_1^2 + y_2^2 = 1, \\ y_1 z_1 z_2 = y_2 \frac{z_1^2 - z_2^2}{2}, \end{cases} \quad \text{and} \quad \begin{cases} y_3^2 + y_4^2 = 1, \\ y_3 z_3 z_4 = y_4 \frac{z_3^2 - z_4^2}{2}. \end{cases}$$

By Theorem 4.6 and Corollary 4.7, we have

$$\sup_{\eta > 0} \inf\{\|x\| : x \in \widehat{D}^* g(z)(y), z \in \mathbb{B}(\bar{z}, \eta), g(z) \in \mathbb{B}(\bar{w}, \eta), \|y\| = 1\}$$

$$= \sup_{\eta > 0} \inf\left\{\|x\| : x \in \widehat{D}^* g(z)(y), z \in \mathbb{B}(\bar{z}, \eta), g(z) \in \mathbb{B}(\bar{w}, \eta), \|y\| = 1, y_1 z_1 z_2 = y_2 \frac{z_1^2 - z_2^2}{2}, y_3 z_3 z_4 = y_4 \frac{z_3^2 - z_4^2}{2}\right\}$$

$$= \sup_{\eta > 0} \inf\left\{\|y\| : x \in \widehat{D}^* f(z)(y), z \in \mathbb{B}(\bar{z}, \eta), g(z) \in \mathbb{B}(\bar{w}, \eta), \|y\| = 1, y_1 z_1 z_2 = y_2 \frac{z_1^2 - z_2^2}{2}, y_3 z_3 z_4 = y_4 \frac{z_3^2 - z_4^2}{2}\right\}$$

$$= \sup_{\eta > 0} \{1\} = 1. \qquad \square$$

Next, we consider a mapping $h : \mathbb{R}^4 \to \mathbb{R}^4$, which is modified the mapping $g : \mathbb{R}^4 \to \mathbb{R}^4$ defined by (4.1). We define $h : \mathbb{R}^4 \to \mathbb{R}^4$, for any $x = (x_1, x_2, x_3, x_4) \in \mathbb{R}^4$, by

$$h((x_1, x_2, x_3, x_4)) = \begin{cases} \left(\frac{x_1^2 - x_2^2}{\|x\|}, \frac{2 x_1 x_2}{\|x\|}, \frac{x_3^2 - x_4^2}{\|x\|}, \frac{2 x_3 x_4}{\|x\|}\right), & \text{if } x \neq \theta, \\ \theta, & \text{if } x = \theta. \end{cases} \qquad (4.8)$$

We see that $h$ is a continuous mapping on $\mathbb{R}^4$. Similar, to the mapping $g : \mathbb{R}^4 \to \mathbb{R}^4$ defined by (4.1), one can study the properties of $h$ at the origin $\theta$. However, in this paper, we only study the Fréchet differentiability and the Mordukhovich differentiability of $h$ on $\mathbb{R}^4 \setminus \{\theta\}$.

**Theorem 4.9.** *Let* $z = (z_1, z_2, z_3, z_4) \in \mathbb{R}^4 \setminus \{\theta\}$. *Then we have that*

(a) *$h$ is Fréchet differentiable at $z$ and*

$$\nabla h(z) = \begin{pmatrix} \frac{(z_1^2 + 3z_2^2 + 2z_3^2 + 2z_4^2) z_1}{\|z\|^3} & \frac{2 z_2^2 z_2}{\|z\|^3} & \frac{-(z_3^2 - z_4^2) z_1}{\|z\|^3} & \frac{-2 z_3 z_4 z_1}{\|z\|^3} \\ \frac{-(3z_1^2 + z_2^2 + 2z_3^2 + 2z_4^2) z_2}{\|z\|^3} & \frac{2 z_1^2 z_1}{\|z\|^3} & \frac{-(z_3^2 - z_4^2) z_2}{\|z\|^3} & \frac{-2 z_3 z_4 z_2}{\|z\|^3} \\ \frac{-(z_1^2 - z_2^2) z_3}{\|z\|^3} & \frac{-2 z_1 z_2 z_3}{\|z\|^3} & \frac{(2z_1^2 + 2z_2^2 + z_3^2 + 3z_4^2) z_3}{\|z\|^3} & \frac{2 z_4^2 z_4}{\|z\|^3} \\ \frac{-(z_1^2 - z_2^2) z_4}{\|z\|^3} & \frac{-2 z_1 z_2 z_4}{\|z\|^3} & \frac{-(2z_1^2 + 2z_2^2 + 3z_3^2 + z_4^2) z_4}{\|z\|^3} & \frac{2 z_3^2 z_3}{\|z\|^3} \end{pmatrix}.$$

(b) *$h$ is Mordukhovich differentiable at $z$ and*

$$\widehat{D}^* h(z) = \begin{pmatrix} \frac{(z_1^2 + 3z_2^2 + 2z_3^2 + 2z_4^2) z_1}{\|z\|^3} & \frac{-(3z_1^2 + z_2^2 + 2z_3^2 + 2z_4^2) z_2}{\|z\|^3} & \frac{-(z_1^2 - z_2^2) z_3}{\|z\|^3} & \frac{-(z_1^2 - z_2^2) z_4}{\|z\|^3} \\ \frac{2 z_2^2 z_2}{\|z\|^3} & \frac{2 z_1^2 z_1}{\|z\|^3} & \frac{-2 z_1 z_2 z_3}{\|z\|^3} & \frac{-2 z_1 z_2 z_4}{\|z\|^3} \\ \frac{-(z_3^2 - z_4^2) z_1}{\|z\|^3} & \frac{-(z_3^2 - z_4^2) z_2}{\|z\|^3} & \frac{(2z_1^2 + 2z_2^2 + z_3^2 + 3z_4^2) z_3}{\|z\|^3} & \frac{-(2z_1^2 + 2z_2^2 + 3z_3^2 + z_4^2) z_4}{\|z\|^3} \\ \frac{-2 z_3 z_4 z_1}{\|z\|^3} & \frac{-2 z_3 z_4 z_2}{\|z\|^3} & \frac{2 z_4^2 z_4}{\|z\|^3} & \frac{2 z_3^2 z_3}{\|z\|^3} \end{pmatrix}.$$

*Proof.* By Lemma 3.1 and Theorem 3.2 in [15], the proof of this theorem is similar to the proof of

Theorem 3.3. It is omitted here. □

From the representation of the Mordukhovich derivative of $h$, it is very difficult and complicated to calculate the covering constant for $h$ by using the similar techniques used in Theorem 4.8. We only consider some special cases, which may be interesting for some readers.

**Proposition 4.10**. *Let* $\bar{z} = (\bar{z}_1, \bar{z}_2, \bar{z}_3, \bar{z}_4) \neq \theta$ *with* $\bar{w} = h(\bar{z})$. *We have*

(a) *If* $\bar{z}_i = 0$, *for some* $i = 1, 2, 3, 4$, *then*

(i) $$\hat{\alpha}(h, \bar{z}, \bar{w}) \leq \frac{2\bar{z}_{3-i}^3}{\|\bar{z}\|^3}, \text{ for } i = 1, 2,$$

(ii) $$\hat{\alpha}(h, \bar{z}, \bar{w}) \leq \frac{2\bar{z}_{7-i}^3}{\|\bar{z}\|^3}, \text{ for } i = 3, 4.$$

*In particular, if* $\bar{z}_1 = \bar{z}_2 = 0$ *or* $\bar{z}_3 = \bar{z}_4 = 0$, *then* $\hat{\alpha}(h, \bar{z}, \bar{w}) = 0$.

(b) *If* $|\bar{z}_1| = |\bar{z}_2| = |\bar{z}_3| = |\bar{z}_4|$, *then*

$$\hat{\alpha}(h, \bar{z}, \bar{w}) \leq \frac{1}{\sqrt{2}}.$$

*Proof.* Proof of (a). Suppose $\bar{z}_1 = 0$. By the continuity of $h$ around $\bar{z} = (\bar{z}_1, \bar{z}_2, \bar{z}_3, \bar{z}_4) \neq \theta$, for any $\eta > 0$, we can have $z \in \mathbb{R}^4 \setminus \{\theta\}$ such that

$$z = (z_1, z_2, z_3, z_4) \in \mathbb{B}(\bar{z}, \eta) \setminus \{\theta\} \text{ satisfying } h(z) \in \mathbb{B}(\bar{w}, \eta) \text{ and } z_1 = 0. \quad (4.9)$$

For any $\eta > 0$, let

$$D(\bar{z}, \eta) = \{z \in \mathbb{R}^4 \setminus \{\theta\} : z \text{ satisfies } (4.9)\}$$

For any given $\eta > 0$, let $z \in \mathbb{R}^4 \setminus \{\theta\}$ satisfying (4.9). Let $x = (x_1, x_2, x_3, x_4), y = (y_1, y_2, y_3, y_4) \in \mathbb{R}^4$. If $x = \widehat{D}^* h(z)(y)$, then, by Theorem 4.9, we have

$$x_1 = y_2 \frac{2z_2^2 z_2}{\|z\|^3}$$

$$x_2 = y_1 \frac{-(z_2^2 + 2z_3^2 + 2z_4^2)z_2}{\|z\|^3} + y_3 \frac{-(z_3^2 - z_4^2)z_2}{\|z\|^3} + y_4 \frac{-2z_3 z_4 z_2}{\|z\|^3}$$

$$x_3 = y_1 \frac{z_2^2 z_3}{\|z\|^3} + y_3 \frac{(2z_2^2 + z_3^2 + 3z_4^2)z_3}{\|z\|^3} + y_4 \frac{2z_4^2 z_4}{\|z\|^3},$$

$$x_4 = y_1 \frac{z_2^2 z_4}{\|z\|^3} + y_3 \frac{-(2z_2^2 + 3z_3^2 + z_4^2)z_4}{\|z\|^3} + y_4 \frac{2z_3^2 z_3}{\|z\|^3}. \quad (4.10)$$

By (4.9) and (4.10), we calculate

$$\hat{\alpha}(h, \bar{z}, \bar{w}) = \sup_{\eta > 0} \inf\{\|x\| : x \in \widehat{D}^* h(z)(y), z \in \mathbb{B}(\bar{z}, \eta) \setminus \{\theta\}, h(z) \in \mathbb{B}(\bar{w}, \eta), \|y\| = 1\}$$

$$\leq \sup_{\eta > 0} \inf\{\|x\| : x \in \widehat{D}^* h(z)(y), z \in (\mathbb{B}(\bar{z}, \eta) \setminus \{\theta\}) \cap D(\bar{z}, \eta), h(z) \in \mathbb{B}(\bar{w}, \eta), \|y\| = 1\}$$

$$\leq \sup_{\eta>0} \inf\{\|x\|: x \in \widehat{D}^*h(z)(y), z \in (\mathbb{B}(\bar{z},\eta)\setminus\{\theta\}) \cap D(\bar{z},\eta), h(z) \in \mathbb{B}(\bar{w},\eta), y_2 = \|y\| = 1\}$$

$$= \sup_{\eta>0} \inf\left\{\frac{2z_2^3}{\|z\|^3}: x \in \widehat{D}^*h(z)(y), z \in (\mathbb{B}(\bar{z},\eta)\setminus\{\theta\}) \cap D(\bar{z},\eta), h(z) \in \mathbb{B}(\bar{w},\eta), y_2 = \|y\| = 1\right\}$$

$$\leq \sup_{\eta>0} \inf\left\{\frac{2\bar{z}_2^3}{\|\bar{z}\|^3}: x \in \widehat{D}^*h(\bar{z})(y), \bar{z} \in (\mathbb{B}(\bar{z},\eta)\setminus\{\theta\}) \cap D(\bar{z},\eta), h(\bar{z}) \in \mathbb{B}(\bar{w},\eta), y_2 = \|y\| = 1\right\}$$

$$= \frac{2\bar{z}_2^3}{\|\bar{z}\|^3}.$$

Rest cases can be similarly proved.

Next, we prove (b). Suppose $|\bar{z}_1| = |\bar{z}_2| = |\bar{z}_3| = |\bar{z}_4|$. Let $x = (x_1, x_2, x_3, x_4) \in \mathbb{R}^4$. We take $y = (0, 1, 0, 0) \in \mathbb{R}^4$ (We may take $y = (0, 0, 0, 1)$). If $x = \widehat{D}^*h(z)(y)$, by Theorem 4.9, we have

$$x_1 = \frac{2z_2^2 z_2}{\|z\|^3}$$

$$x_2 = \frac{2z_1^2 z_1}{\|z\|^3}$$

$$x_3 = \frac{-2z_1 z_2 z_3}{\|z\|^3},$$

$$x_4 = \frac{-2z_1 z_2 z_4}{\|z\|^3}. \qquad (4.11)$$

By (4.11), we calculate

$$\hat{\alpha}(h, \bar{z}, \bar{w}) = \sup_{\eta>0} \inf\{\|x\|: x \in \widehat{D}^*h(z)(y), z \in \mathbb{B}(\bar{z},\eta)\setminus\{\theta\}, h(z) \in \mathbb{B}(\bar{w},\eta), \|y\| = 1\}$$

$$\leq \sup_{\eta>0} \inf\{\|x\|: x \in \widehat{D}^*h(z)(y), z \in (\mathbb{B}(\bar{z},\eta)\setminus\{\theta\}) \cap D(\bar{z},\eta), h(z) \in \mathbb{B}(\bar{w},\eta), y_2 = \|y\| = 1\}$$

$$\leq \sup_{\eta>0} \inf\left\{\sqrt{\frac{16\bar{z}_1^6}{(2|\bar{z}_1|)^6}}: x \in \widehat{D}^*h(\bar{z})(y), \bar{z} \in (\mathbb{B}(\bar{z},\eta)\setminus\{\theta\}) \cap D(\bar{z},\eta), h(\bar{z}) \in \mathbb{B}(\bar{w},\eta), y_2 = \|y\| = 1\right\}$$

$$= \frac{1}{\sqrt{2}}. \qquad \square$$

### 5. Some Applications

One of the important applications of the Mordukhovich derivatives is to define and to calculate the covering constants for some set-valued mappings in Banach spaces. By the concept of covering constants, in [1], the authors proved the Arutyunov Mordukhovich and Zhukovskiy Parameterized Coincidence Point Theorem (It is simply named as AMZ Theorem). The AMZ Theorem provides a very powerful tool in set-valued and variational analysis (see [1, 11−15]). In this section, we use the results $\hat{\alpha}(f, \bar{z}, \bar{w}) = 1$ and apply the AMZ Theorem to solve some parameterized equations in $\mathbb{R}^2$ related to the single-valued mapping $f: \mathbb{R}^2 \to \mathbb{R}^2$ defined by (1.1).

**(AMZ Theorem)** *Let the Banach spaces X and Y in be Asplund and let P be a topological space. Let $F: X \rightrightarrows Y$ and $G(\cdot, \cdot): X \times P \rightrightarrows Y$ be set-valued mappings. Let $\bar{x} \in X$ and $\bar{y} \in Y$ with $\bar{y} \in F(\bar{x})$. Suppose*

*that the following conditions are satisfied:*

**(A1)** *The multifunction $F: X \rightrightarrows Y$ is closed around $(\bar{x}, \bar{y})$.*

**(A2)** *There are neighborhoods $U \subset X$ of $\bar{x}$, $V \subset Y$ of $\bar{y}$, and $O$ of $\bar{p} \in P$ as well as a number $\beta \geq 0$ such that the multifunction $G(\cdot, p): X \rightrightarrows Y$ is Lipschitz-like on $U$ relative to $V$ for each $p \in O$ with the uniform modulus $\beta$, while the multifunction $p \to G(\bar{x}, p)$ is lower/inner semicontinuous at $\bar{p}$.*

**(A3)** *The Lipschitzian modulus $\beta$ of $G(\cdot, p)$ is chosen as $\beta < \hat{\alpha}(F, \bar{x}, \bar{y})$, where $\hat{\alpha}(F, \bar{x}, \bar{y})$ is the covering constant of $F$ around $(\bar{x}, \bar{y})$ taken from (2.5).*

*Then for each $\alpha > 0$ with $\beta < \alpha < \hat{\alpha}(F, \bar{x}, \bar{y})$, there exist a neighborhood $W \subset P$ of $\bar{p}$ and a single-valued mapping $\sigma: W \to X$ such that whenever $p \in W$ we have*

$$F(\sigma(p)) \cap G(\sigma(p), p) \neq \emptyset \quad \text{and} \quad \|\sigma(p) - \bar{x}\|_X \leq \frac{\text{dist}(\bar{y}, G(\bar{x}, p))}{\alpha - \beta}.$$

**Theorem 5.1.** *Let $C$ be a topological space. Let $f: \mathbb{R}^2 \to \mathbb{R}^2$ be a single-valued mapping defined by (1.1). and let $h(\cdot, \cdot) = (h_1(\cdot, \cdot), h_2(\cdot, \cdot)): \mathbb{R}^2 \times C \to \mathbb{R}^2$ be single-valued mappings. Let $\omega = (\omega_1, \omega_2): C \to \mathbb{R}^2$ be a single-valued lower semicontinuous mapping. Let $\bar{x}$ and $\bar{y} \in \mathbb{R}^2 \setminus \{\theta\}$ with $\bar{y} = f(\bar{x})$. Suppose that the following conditions are satisfied:*

*There are neighborhoods $U \subset \mathbb{R}^2$ of $\bar{x}$, $V \subset \mathbb{R}^2$ of $\bar{y}$ and $O$ of $\bar{s} \in C$ as well as a number $\beta \in [0, 1)$ such that the mapping $h(\cdot, s): \mathbb{R}^2 \to \mathbb{R}^2$ satisfies the Lipschitz condition on $U$ relative to $V$ for each $s \in O$ with the uniform modulus $\beta$, while the mapping $s \to h(\bar{x}, s)$ is lower semicontinuous at $\bar{s}$.*

*Then for each $\alpha > 0$ with $\beta < \alpha < 1$, there exist a neighborhood $W \subset C$ of $\bar{s}$ and a single-valued mapping $\sigma = (\sigma_1, \sigma_2): W \to \mathbb{R}^2$ such that whenever $s \in W$ we have*

$$\begin{cases} \frac{\sigma_1(s)^2 - \sigma_2(s)^2}{\sqrt{\sigma_1(s)^2 + \sigma_2(s)^2}} = h_1(\sigma_1(s), \sigma_2(s)) + \omega_1(s), \\ \frac{2\sigma_1(s)\sigma_2(s)}{\sqrt{\sigma_1(s)^2 + \sigma_2(s)^2}} = h_2(\sigma_1(s), \sigma_2(s)) + \omega_2(s), \end{cases}$$

*and*

$$\|(\sigma_1(s), \sigma_2(s)) - (\bar{x}_1, \bar{x}_2)\| \leq \frac{\left\| h((\bar{x}_1, \bar{x}_2), s) + (\omega_1(s), \omega_2(s)) - \left( \frac{\bar{x}_1^2 - \bar{x}_2^2}{\sqrt{\bar{x}_1^2 + \bar{x}_2^2}}, \frac{2\bar{x}_1 \bar{x}_2}{\sqrt{\bar{x}_1^2 + \bar{x}_2^2}} \right) \right\|}{\alpha - \beta}. \quad (4.1)$$

*Proof.* In AMZ Theorem, we particularly take

$$F(x) = f(x) \text{ and } G(x, s) = h(x, s) + \omega(s), \text{ for all } x \in \mathbb{R}^2, s \in C. \quad (5.2)$$

By Theorem 3.5, we have that

$$\hat{\alpha}(f, \bar{z}, f(\bar{z})) = 1, \text{ for any } \bar{z} \in \mathbb{R}^2.$$

As single-valued mappings, we can verify that the mappings, which are defined in (5.2), $F: \mathbb{R}^2 \to \mathbb{R}^2$ and $G(\cdot, \cdot): \mathbb{R}^2 \times C \to \mathbb{R}^2$ satisfy all conditions in the AMZ Theorem, from which the conclusions (4.1) are immediately proved. □

**Acknowledgments.** The author is very grateful to Professor Boris S. Mordukhovich, Professor Simeon Reich, Professor Christiane Tammer and Professor J. C. Yao for their kind communications, valuable suggestions and enthusiasm encouragements in the development stage of this paper.